\documentclass[12pt]{article}
\usepackage[english]{babel}
\usepackage{amsmath,amsthm,amsfonts,amssymb,epsfig,vector, simplemargins,mathrsfs}
\numberwithin{equation}{section}

\setallmargins{1in}

\newtheorem{thm}{Theorem}[section]
\newtheorem*{thmdim}{Theorem \ref{thm:dim}}
\newtheorem*{thmsphere}{Theorem \ref{thm:sphere}}
\newtheorem{lem}[thm]{Lemma}
\newtheorem{prop}[thm]{Proposition}
\newtheorem{cor}[thm]{Corollary}
\newtheorem{df}[thm]{Definition}

\newtheorem{conj}{Conjecture}

\newcommand{\M}{\mathcal{M}}
\newcommand{\B}{\mathcal{B}}

\newcommand{\hilbert}{\mathcal{H}}
\newcommand{\G}{\mathcal{G}}
\newcommand{\prob}{\mathbb{P}}
\newcommand{\E}{\mathbb{E}}
\newcommand{\R}{\mathbb{R}}
\newcommand{\Z}{\mathbb{Z}}
\newcommand{\N}{\mathbb{N}}
\newcommand{\Q}{\mathbb{Q}}

\newcommand{\distrib}{\overset{\mbox{law}}{=}}

\newcommand{\ssigma}{\vec{\sigma}}
\newcommand{\bbeta}{\vec{\beta}}

\begin{document}
\author{{\Large Louis-Pierre Arguin}
\thanks{L.-P. Arguin held a postdoctoral position at Courant Institute during this work. He was supported by the NSF grant DMS-0604869 and partially by the Hausdorff Center for Mathematics, Bonn.} 
\\[1ex] 
\small D\'epartement de Math\'ematiques et Statistique, Universit\'e de Montr\'eal\\
\small Montr\'eal, Qu\'ebec, H3T 1J4, Canada
\\[2ex] 
{\Large Sourav Chatterjee}
\thanks{Sourav Chatterjee's research was partially supported by NSF grants DMS-0707054 and  DMS-1005312,   and a Sloan Research Fellowship}
\\[1ex] 
\small Courant Institute of Mathematical Sciences\\
\small New York University, New York, 10012, USA.}

\date{April 16, 2012}

\title{Random Overlap Structures: \\Properties and Applications to Spin Glasses}

\maketitle

\begin{abstract}
Random Overlap Structures (ROSt's) are random elements on the space of probability measures on the unit ball of a Hilbert space, where two measures are identified if they differ by an isometry. In spin glasses, they arise as natural limits of Gibbs measures under the appropriate algebra of functions. We prove that the so called `cavity mapping' on the space of ROSt's is continuous, leading to a proof of the stochastic stability conjecture for the limiting Gibbs measures of a large class of spin glass models. Similar arguments yield the proofs of a number of other properties of ROSt's that may be useful in future attempts at proving the ultrametricity conjecture. 
Lastly, assuming that the ultrametricity conjecture holds, 
the setup yields a constructive proof of the Parisi formula for the free energy of the Sherrington-Kirkpatrick model by making rigorous a heuristic of Aizenman, Sims and Starr.
\end{abstract}

\section{Introduction}
This paper develops some connections between random probability measures on Hilbert spaces and the Gibbs measures of spin glasses with Gaussian disorder. Consider for $N\in \N$ random processes of the form
$$
H_N:=(H_{N}(\sigma), \sigma\in \{-1,+1\}^N)\ ,
$$
where $H_N(\sigma)$ is a centered Gaussian variable. Denote the law of $H_N$ by $P$ and the corresponding expectation by $E$.
Suppose that the covariances between variables are of the form
$$
E \ H_{N}(\sigma) H_{N}(\sigma')=N r(\sigma,\sigma')$$ 
for some positive definite symmetric form $r$ on $\{-1,+1\}^N$ with $r(\sigma,\sigma)=1$ for all $\sigma$.
We say that $H_N$ is the {\it Hamiltonian of a Gaussian spin glass}.
An important example is the Sherrington-Kirkpatrick (SK) model for which 
$$
 r(\sigma,\sigma')=\left(\frac{1}{N}\sum_{i=1}^N\sigma_i\sigma_i'\right)^2\ .
$$
The definition also includes the Edwards-Anderson (EA) model. 
This model can be defined for example in a box of $\Z^d$ with $N$ vertices and edge set $\mathcal{E}_N$.
The form for the EA model is then 
$$
r(\sigma,\sigma')=\frac{1}{|\mathcal{E}_N|}\sum_{\{i,j\}\in \mathcal{E}_N}\sigma_i\sigma_j\sigma_i'\sigma_j' \ .
$$
The Gibbs measure corresponding to the Hamiltonian  $H_N$ at inverse temperature $\beta$ is defined by
\begin{equation} \label{eqn: Gibbs measure}
\mathcal{G}_{\beta,N}(\sigma)=\frac{\exp \beta H_N(\sigma)}{Z_N(\beta)}
\end{equation}
where $Z_N(\beta)=\sum_{\sigma} \exp \beta H_N(\sigma)$.  
A fundamental problem in the theory of spin glasses is to describe
the limits of $\mathcal{G}_{\beta,N}$ as $N$ tends to infinity. 
This information is useful in particular to understand the extreme value statistics
of the Gaussian process $H_N$.
In the past ten years, the rigorous study of spin glasses has made important progresses 
in the understanding of the models in the limit $N\to\infty$. One of the major achievements was the proof of the Parisi formula
by Guerra and Talagrand \cite{guerra_bound,tal_proof} for the limiting free energy (that is, $\lim_{N\to\infty}\frac{1}{N}E\log Z_N(\beta)$) of the SK model.
The formula suggests that the limiting Gibbs measure $\mathcal{G}_{\beta}$ 
of the SK model (in a suitable sense) has a support that is {\it hierarchical} or {\it ultrametric}, that is with $P$-probability one
$$
\mathcal{G}_{\beta}^{\otimes 3}\Big\{(\sigma,\sigma',\sigma''):r(\sigma,\sigma')\geq \min\left\{r(\sigma,\sigma''), r(\sigma',\sigma'')\right\} \Big\}=1\ .
$$
In fact, it is expected that the ultrametricity of the Gibbs measure is universal to a certain extent within spin glasses. 
This hypothesis is at the core of the description of spin glasses developed by physicists, and often referred to as Parisi theory \cite{MVP}.
In spite of important advances in the understanding of the structure of the Gibbs measure \cite{panchenko_GG, argaiz},
the problem of rigorously characterizing the limiting Gibbs measures of spin glasses remains open. 

This paper takes the approach introduced by Aizenman, Sims and Starr~\cite{ASS} and extended in \cite{lp_remark} to study the limiting Gibbs measures of spin glasses. 
The idea goes as follows. Given an infinite-dimensional separable Hilbert space $\hilbert$ with inner product ``$\cdot$'', there exists an embedding of the hypercube $\{-1,1\}^N$ in $\hilbert$ 
that preserves the form $r$, that is, if $\sigma\mapsto v(\sigma)\in\hilbert$ then
$$
r(\sigma,\sigma')=v(\sigma)\cdot v(\sigma')\ .
$$ 
For example, for the SK model with $N$ sites, the Hilbert space $\hilbert$ can be taken to be the $L^2$ space of some abstract probability space $(\Omega, \mathcal{F}, \prob)$, and the map $v$ defined as
\[
v(\sigma) = \frac{1}{N}\sum_{i,j=1}^N g_{ij} \sigma_i \sigma_j,
\]
where $g_{ij}:\Omega \rightarrow \mathbb{R}$ are independent standard Gaussian random variables under $\prob$. 

Under such an embedding, the Gibbs measure is naturally
sent to the random probability measure on $\hilbert$ assigning weight $\G_{\beta,N}(\sigma)$ to the corresponding vector $v(\sigma)$.
Of course, there is more than one such isometric embedding.
Since the products $v\cdot v'$ determine the vectors up to a choice of basis for $\hilbert$, 
it is necessary to identify the probability measures on $\hilbert$ that differ by an isometry for
the image measure of  $\mathcal{G}_{\beta,N}$ to be well-defined.
It turns out that the proper way to study such measures is through the notion of  a 
{\it Random Overlap Structure} or ROSt. 
A ROSt is a $\N\times\N$ random covariance matrix whose law is invariant under simultaneous finite permutations of 
rows and columns. In the spin glass setting, the ROSt is the matrix $Q^{\beta,N}=\{r(\sigma^k,\sigma^l)\}_{k,l\in\N}$ where $(\sigma^k,k\in\N)$ 
are sampled independently from $\G_{\beta,N}$.
The space of ROSt's is compact in the appropriate topology 
and it is in correspondence, up to isometry, with the space of random probability measures on $\hilbert$
(see Section \ref{section: rost}).
Thus one way to study the limiting Gibbs measures of spin glasses is 
to describe the limit points of the sequence $(Q^{\beta,N})$ in the space of ROSt's. 
(Another promising approach is through the Ghirlanda-Guerra (GG) identities \cite{guerra_GG, panchenko_GG, tal_pure}. 
Although the GG identity approach is not the focus of this article, it seems to be closely connected to the ROSt approach at some level.)

The main result of this paper is to establish that the limiting Gibbs measures (in the ROSt sense) of any Gaussian spin glass as defined above
is {\it stochastically stable}.
What this means is that the ROSt's for Gaussian spin glasses are distributionally invariant under a built-in stochastic mapping that we call the {\it cavity mapping}. (The mapping is well-known, but the nomenclature is ours. Precise definition is given below.)  
The transformation derives its name and structure from  the cavity method introduced in \cite{mvp_cavity}
(see also \cite{BS} for the definition of the mapping for Ruelle Probability Cascades).
It has the features of a basic stochastic object akin to the mappings studied by Kahane in the context of multiplicative chaos \cite{kahane}.
Similar mappings have been considered earlier  in the setting of competing particle systems (see e.g.~\cite{ruzmaiz,argaiz}).

It was suggested by the work of several authors that the Gibbs measures of spin glasses must be stochastically stable \cite{AC, ASS, ruzmaiz, argaiz}. Such a statement, however, can be made precise only in the infinite-dimensional setting. One missing ingredient to make the assertion rigorous at the level of ROSt's was a proof of the continuity of the cavity mapping.
This is done in Theorem \ref{thm:continuity}. Combining this with some further ingredients, it follows that the limiting Gibbs measure of a Gaussian spin glass is stochastically
stable at any $\beta$ where the free energy is differentiable (by convexity, this holds for almost all~$\beta$).

The characterization of stochastically stable ROSt's is a challenging problem.
In view of the main result of this paper, it is also a rewarding one since it provides a way to establish
universal properties of the Gibbs measures of spin glasses, in particular the alleged ultrametricity.
We formulate below a precise conjecture, which may be called the Ultrametricity Conjecture for ROSt's: a ROSt that is stochastically stable
(in a strong sense to be defined) must have ultrametric support.
Some partial results towards the resolution of this conjecture are given in Section \ref{section:properties}.
The most important one might be Theorem \ref{thm:dim} which states that the support of a stochastically stable ROSt is either one vector or lies
in an infinite-dimensional subspace almost surely.

The framework of ROSt's has other interesting consequences for spin glasses. In particular, 
assuming that the ultrametricity conjecture for ROSt's is true, we give a constructive proof of the Parisi formula for the free energy of the Sherrington-Kirkpatrick model. 
This makes rigorous the approach suggested by Aizenman, Sims and Starr \cite{ASS} to prove the Parisi formula.

We remark that results similar to ours have been proved in parallel by Panchenko \cite{panchenko_spins}. 
The work focuses on mean-field spin glass models where the law of the Hamiltonian is invariant under permutation of spins,
such as the SK model. In this case, a representation theorem for exchangeable arrays (the Aldous-Hoover theorem)
can be used. In the ROSt setting, we appeal to the Dovbysh-Sudakov theorem (Theorem \ref{thm:DS} below).
A continuity result is proved and a representation for the Parisi formula, similar to the one obtained here for ROSt's, is found

The paper is organized as follows. The framework of ROSt's as well as the precise statements of the main results are given in Section \ref{section: rost}.
The applications of the framework to spin glasses are  presented in Section \ref{section: spin glass}. 
The continuity of the cavity mapping is proved in Section \ref{section:cavity}. Section \ref{section: parisi formula}
derives the Parisi formula from the Ultrametricity Conjecture. Finally, some interesting properties of stochastically stable ROSt's are proved in Section \ref{section:properties}.

{\bf Acknowledgments.} We thank the referee for insightful remarks and observations that led to major improvements of the first version of this paper.

\subsection{Random Overlap Structures and Probability Measures on Hilbert Spaces}
\label{section: rost}
This section begins with the definition of random overlap structures using exchangeable covariance matrices. The definition is then explained in terms of probability measures on Hilbert spaces using
the representation theorem of Dovbysh and Sudakov \cite{DS}. This setup was previously described in \cite{lp_remark}. The section ends with the definition of the cavity mapping and the statement of several related results. 

Consider the space of positive semi-definite symmetric $\N\times\N$ matrices, or covariance matrices, with $1$ on the diagonal.
This is a compact separable metric space when considered as a closed subset of $[-1,1]^{\N\times\N}$ equipped with the product topology. 
The Borel probability measures on this space form also a compact, separable and metric space when equipped with the weak-$*$ topology 
generated by the continuous functions on the entries. An element $\prob$ of this space is the law of some random covariance matrix $Q=\{q_{ij}\}$.
The continuous functionals in the topology considered can be approximated by linear combinations of monomials of the form
\begin{equation}
\prob \mapsto \E\prod_{1\leq i<j\leq s}q_{ij}^{k_{ij}}\ ,
\label{eqn:fct}
\end{equation}
for $s\in\N$ and a collection $k_{ij}\in\N$.
A random covariance matrix $Q$ is said to be {\it weakly exchangeable} if for any permutation matrix $\tau$ that fixes all but a finite number of elements
$$ \tau \ Q \ \tau^{-1}\ \distrib Q\ ,$$
where $\distrib$ denotes equality in law.
Since this symmetry is preserved under convergence defined by the functions \eqref{eqn:fct},
the space of distributions of weakly exchangeable random covariance matrices with $1$ on the diagonal is a compact 
convex set. Throughout the paper, we refer to this subset of laws as the {\it space of ROSt's}.
\begin{df}
A random overlap structure, or ROSt, is a weakly exchangeable $\N\times \N$ random covariance matrix with $1$ on the diagonal.
The law of a ROSt will be usually denoted by $\prob$ and integration with respect to $\prob$ by $\E$.
\end{df}
Since the law of a ROSt is determined by the countable set of functions of the form \eqref{eqn:fct},
it is convenient to see these as the coordinates of the ROSt.

The relation between ROSt's and probability measures on Hilbert spaces is provided by the 
representation theorem of Dovbysh and Sudakov \cite{DS} (see \cite{panchenko_DS} for a readable proof of this result).
Let $\B$ be the unit ball of an infinite-dimensional separable Hilbert space $\hilbert$ (fixed throughout this paper).
The inner product on $\hilbert$ will be denoted by $\cdot$ and a generic element of $\mathcal{B}$ by~$v$. 
The set of Borel probability measures on $\B$ is denoted by $\M(\B)$.
Note that $\M(\B)$ is the same whether $\B$ is equipped with the weak topology or the norm topology.
From any element of $\M(\B)$, one can construct a ROSt by taking
iid vectors sampled from this element and considering the Gram matrix of these vectors. 
The Dovbysh-Sudakov theorem states that any weakly exchangeable random covariance matrix can be 
constructed this way if one considers random elements of $\M(\B)$. 
In the case where the diagonal entries of the matrices are all $1$, the Dovbysh-Sudakov theorem reads: 
\begin{thm}
Let $Q=\{q_{ij}\}$ be a ROSt. There exists a random element $\mu$ of $\M(\B)$ such that 
conditionally on $\mu$,
$$
(q_{ij})_{i\ne j}\distrib (v^i\cdot v^j)_{i\ne j}, 
$$
where $( v^i,i\in\N)$ are iid $\mu$-distributed vectors in $\B$.
\label{thm:DS}
\end{thm}
We refer to the measure $\mu$ as the {\it sampling measure} of the ROSt $Q$ in this article. 
It plays a role analogous to the empirical measure in de Finetti's theorem.
The law of $Q$ does not determine the sampling measure uniquely. However, there is a one-to-one correspondence between ROSt's and random elements of $\M(\B)$ if
the elements of $\M(\B)$ that differ by an isometry are identified. The reader is referred to \cite{panchenko_DS} for a proof.
\begin{prop}
\label{prop: sampling}
Let $Q$ and $Q'$ be two ROSt's with sampling measures $\mu$ and $\mu'$. Then $Q\distrib Q'$ if and only if
there an isometry $T$ of $\B$ (possibly dependent on the realization of $\mu$)
such that $\mu(A) \distrib \mu'(T^{-1} A)$ for any $A\subset \B$ Borel measurable. 
\end{prop}
We write $\vec{v}=( v^1,..., v^s)$ for an element of $\B^s$ and 
$\mu^{\otimes s}$ for the product measure of $s$ copies of $\mu$. Here $s$ is an arbitrary positive integer. 
Each sampled vector is called a {\it replica}.
Throughout the paper, expectation of a function $F:\B^s\to\R$ with respect to $\mu^{\otimes s}$ 
will often be written $\mu^{\otimes s}(F)$ for short.
In this notation, the continuous functionals \eqref{eqn:fct} for ROSt's are all of the form
\begin{equation}
\prob\mapsto \E\left[\mu^{\otimes s}\left(F( \vec{v})\right)\right]\ .
\label{eqn:h}
\end{equation}
where $F:\B^s\to \R$ is a {\it continuous function on $s$ replicas}, that is, a continuous function that 
depends only on the inner product between $s$ distinct vectors:
$$F( \vec{v})=F( v^k\cdot v^l; 1\leq k<l\leq s)\ .$$
The central topic of interest in this article is a class of ROSt's that possess an invariance property called {\it stochastic stability}. This is defined as follows. 

Let $\prob$ be a ROSt. Let $P_l$ be the law of a collection of independent standard Gaussian random variables $(l_i,i\in\N)$, defined on the same probability space as $\prob$ and independent of it. Let $E_l$ denote integration with respect to $P_l$. For each $v\in\B$, we consider the Gaussian random variable $l(v)=\sum_{i\in\N}l_i \  v\cdot e_i$, where $(e_i,i\in\N)$ is a fixed orthonormal basis of $\hilbert$. 
The covariance of $l(v)$ and $l(v')$ is given by
$$
E_l[l( v)l( v')]= v\cdot v'\ .
$$
The Gaussian process $l=(l( v), v\in\mathcal{B})$ indexed by $\B$ is called the {\it cavity field} in this article. 

Note that the function $ v\mapsto l( v)=\sum_{i\in\N}l_i \  v\cdot e_i$ on $\B$ is measurable $P_l$-almost surely
being the pointwise limit of weakly continuous functions.
Moreover, because $\E E_l \left[\mu(l( v)^2)\right]\leq 1$ by Fubini's theorem, for any law of ROSt $\prob$ there exists a set of realizations of $l$ of $P_l$-probability one such that $l( v)$ is well-defined for $\mu$-almost all $ v$ and for $\prob$-almost all $\mu$.

We now define the {\it cavity mapping} from the space of ROSt's to itself. 
The mapping is defined on the sampling measure.
Let $\lambda>0$ be a parameter, 
$\prob$ be the law of a ROSt with sampling measure $\mu$, and $l$ be a cavity field independent of $\prob$, we take:
\begin{equation}\label{eqn:mapping}
\Phi_{\lambda, l}:\mu \mapsto \Phi_{\lambda, l}\mu:=\frac{e^{\lambda l( v)-\frac{\lambda^2}{2}\| v\|^2}\ \mu(d v) }{\mu(e^{\lambda l( v)-\frac{\lambda^2}{2}\| v\|^2})}\ .
\end{equation}
The image ROSt is constructed using $\Phi_{\lambda, l}\mu$ as sampling measure. Its law is induced by $\prob\times P_l$.
We stress that the factor $e^{\frac{\lambda^2}{2}\| v\|^2}$ is the expectation of $e^{\lambda l( v)}$, and its presence turns
out to be crucial for the continuity of the mapping.
The cavity mapping can also be defined in terms of the continuous functions \eqref{eqn:h} since they determine the law of the ROSt:
$$
\E[\mu^{\otimes s}(F( \vec{v}))]\mapsto \E[(\Phi_{\lambda, l}\mu)^{\otimes s}(F( \vec{v}))]=\E E_l\left[ \frac{\mu^{\otimes s}(F( \vec{v})e^{\lambda l( v^1)-\frac{\lambda^2}{2}\| v^1\|^2}...e^{\lambda l( v^s)-\frac{\lambda^2}{2}\| v^s\|^2})}{\mu^{\otimes s}(e^{\lambda l( v^1)-\frac{\lambda^2}{2}\| v^1\|^2}...e^{\lambda l( v^s)-\frac{\lambda^2}{2}\| v^s\|^2})}\right]\ .
$$
It will be a consequence of Theorem \ref{thm:continuity} that the cavity mapping is well-defined, 
that is, the image ROSt is the same if we apply \eqref{eqn:mapping} to two random elements of $\M(\B)$ that are sampling measures
of the same ROSt.

We remark that this mapping can be generalized in at least two ways.
First, for any function $\psi\in\mathcal{C}^1(\R)$ with bounded derivative, we can take the change of 
density $E_z[e^{\psi(l( v)+z\sqrt{1-\| v\|^2})}]$, where $z$ is a standard Gaussian with expectation $E_z$.
Second, one could consider the cavity field 
$l_c$ with $E_l[l_c( v)l_c( v')]=c( v\cdot v')$ for some continuous function $c$ satisfying
\begin{equation}
\label{eqn:c}
\begin{aligned}
c(1)=1 \hspace{0.4cm} \text{$ ( v, v')\mapsto c( v\cdot v')$ is positive definite}\ .
\end{aligned}
\end{equation}
We will often write $l_p$ corresponding to the
choice $c( v\cdot v')=( v\cdot v')^p$ for $p\in\N\cup\{0\}$.
Equation \eqref{eqn:mapping} is the particular case $\psi(x)=\lambda x$ and $c(x)=x$.
The generalized mapping takes the form 
\begin{equation}
\label{eqn:mapping2}
 \Phi_{\psi,l_c}:\mu\mapsto \Phi_{\psi,l_c}  \mu:= \frac{E_z[e^{\psi(l_c( v)+z\sqrt{1-c(\| v\|^2)})}]\ \mu(d v) }{\mu(E_z[e^{\psi(l_c( v)+z\sqrt{1-c(\| v\|^2)})}])}\ .
\end{equation}

The first result of this paper, from which all others essentially stem, is proved in Section \ref{section:cavity}.
\begin{thm}\label{thm:continuity}
The mapping from the space of ROSt's to itself defined by \eqref{eqn:mapping2} is continuous.
\end{thm}
We are interested in ROSt's whose law is invariant under the cavity mapping \eqref{eqn:mapping2}.
This is equivalent to say by Proposition \ref{prop: sampling} that the
sampling measure $\mu$ has the same law as $\Phi_{\psi,l_c}\mu$, up to an isometry of $\B$.
One direct consequence of the continuity is that the set of invariant laws under a family of cavity mappings
is compact in the space of laws of ROSt's (see Corollary \ref{cor:compact}), a simple fact that will turn out to be useful,
for example in the proof of the Parisi formula (see Corollary \ref{cor: beta positive}).
Our focus is on ROSt's whose law is invariant under the cavity mapping for linear functions $\psi$.
\begin{df}
\label{df: SS}
A ROSt of law $\prob$ is called stochastically stable if, for any $\lambda>0$, its law is invariant under the cavity mapping \eqref{eqn:mapping}.
In other words, for any $s\in\N$ and any continuous function $F$ on $s$ replicas
$$
\E[\mu^{\otimes s}(F( \vec{v}))]= \E[(\Phi_{\lambda, l}\mu)^{\otimes s}(F( \vec{v}))]\ .
$$ 
It is called stochastically stable for the cavity field $l_c$ if the above holds with $l$ replaced
by the cavity field $l_c$ with covariance $E_l[l_c( v)l_c( v')]=c( v\cdot v')$.
\end{df}
As we shall explain in the next section, continuity provides the missing link between
the stochastic stability as defined here and stochastic stability of the Gibbs measure of spin glasses
as studied in \cite{AC,CG,tal_pure}. This link provides a statement of the ultrametricity conjecture on the characterization of stochastically stable ROSt's.
\begin{conj}[Ultrametricity Conjecture]
\label{conj:ultra}
If a ROSt of law $\prob$ is stochastically stable for the cavity fields $l_p$ (defined below \eqref{eqn:mapping2}) for $p=1$ and infinitely many $p$, 
then the support of its sampling measure $\mu$ is ultrametric $\prob$-a.s., that is:
$$
 v\cdot v'\geq \min\left\{ v\cdot v'',  v'\cdot v''\right\} \text{ \ for $\mu^{\otimes 3}$-almost all $ v, v', v''$.}
$$
Moreover, the law of the sampling measure is a convex combination of Ruelle Probability Cascades.
\end{conj}
Ruelle Probability Cascades (RPC's) form a compact subset of ROSt's parametrized by the 
right-continuous, increasing functions on $[0,1]$ with values $0$ at $0$ and $1$ at $1$ (we refer to \cite{lp_remark} for a proof of the compactness). They have ultrametric support and
they are known to be stochastically stable for any cavity field $l_p$, $p\in\N$.
We shall not need their precise definition here. The reader is referred to \cite{ruelle, BS, argaiz} for more details.
It is very likely that the hypotheses of the conjecture can be relaxed. However, it does not hold
if stochastic stability for a single cavity field is assumed. A simple counterexample is given at the end of Section \ref{section:properties}.
The conjecture was proved in the case where the support of $\mu^{\otimes 2}$ is finite in \cite{argaiz} (see also \cite{miller} for an extension
of that proof). In this case,
$\mu$ is supported on countable vectors and the conjecture can be formulated in the language of competing particle systems.
A similar statement of the ultrametricity conjecture exists where the assumption of stability is replaced with the extended Ghirlanda-Guerra identities:
\begin{df}
A ROSt with sampling measure $\mu$ is said to satisfy the Ghirlanda-Guerra (GG) Identities if for any $s\in\N$ and any continuous function $F$
on $s$ replicas we have
\[
\begin{aligned}
&\E[\mu^{\otimes s+1}( v^1\cdot  v^{s+1} F( \vec{v}))]\\
&\hspace{2cm}=\frac{1}{s}\E[\mu^{\otimes 2}( v^1\cdot v^2)]\ \E[\mu^{\otimes s}( F( \vec{v}))]+
\frac{1}{s}\sum_{l=2}^s\E[\mu^{\otimes s}( v^1\cdot v^l\ F( \vec{v}))]\ .
\end{aligned}
\]
It is said to satisfy the extended Ghirlanda-Guerra identities if the above holds when the products $ v^i\cdot v^j$ in the expectations is replaced 
by $( v^i\cdot v^j)^p$ for any $p\in\N$.  
\end{df}
A version of the ultrametricity conjecture assuming finiteness of the support of the overlap distribution and the extended Ghirlanda-Guerra identities was proved by Panchenko in the case where $\mu^{\otimes 2}$ has finite support  \cite{panchenko_GG}. 
(The assumptions are slightly different from the ones used in the proof of the conjecture using stochastic stability in \cite{argaiz} where stability for infinitely many $p$'s in $\N$, but not all $p$, is needed.
However, it is possible that the two sets of hypotheses are equivalent in the special case where the support is finite.)
We also stress that ultrametricity has been proved for generalized versions of GREM spin glass models by Bolthausen and Kistler
\cite{BK1,BK2}.

The continuity of the cavity mapping has interesting consequences that we discuss in Section \ref{section:properties}. 
For one, it provides an ergodic decomposition of stochastically stable ROSt's. 
The Ultrametricity Conjecture can then be restated by saying that the extremes of the considered set of 
stochastically stable ROSt's are exactly the RPC's.
The decomposition is useful to prove properties of the support of the sampling measure. 
We stress that invariance under a single cavity mapping is sufficient here.
\begin{thm}
\label{thm:dim}
Let $\prob$ be the law of a ROSt that is invariant under the cavity mapping $\Phi_{\lambda,l}$
for a given $\lambda>0$. Let $\mu$ be its sampling measure.
Then $\mu$ is supported on a single vector or on an infinite-dimensional subset of $\B$.
\end{thm}
If invariance under more fields is assumed, we can prove more.

\begin{thm}
\label{thm:sphere}
Let $\prob$ be the law of a ROSt that is invariant under the cavity mapping $\Phi_{\lambda,l_p}$
for a given $\lambda>0$ and for infinitely many $p\in\N$.
Let $\mu$ be its sampling measure and suppose that $\mu^{\otimes 2}\{(v,v')\in\B^2:  |v\cdot v'|=r^2_{\max}(\mu)\}>0$ 
on a set of positive $\prob$-probability, where $r_{\max}(\mu)=\sup\{0\leq r \leq 1: \mu\{\|v\|\geq r\}>0\}$.
Then,
$$
\mu\{v\in\B: \|v\|=r_{\max}(\mu)\}=1\ \text{ $\prob$-a.s.}
$$
In other words,  $\mu$ is almost surely supported on a sphere.
\end{thm}
This last property does not hold in general if stability is assumed for a single cavity field. A counter-example is given in Section \ref{section:properties}.

Properties of the sampling measure for ROSt's satisfying the extended Ghirlanda-Guerra identities have been obtained in \cite{panchenko_GG}. 
The apparent similarities of the properties of the ROSt's satisfying the Ghirlanda-Guerra identities
and stochastic stability motivate the following question: are the Ghirlanda-Guerra identities and stochastic stability two
representations of the same property ? 
Since the set of stochastic stable ROSt's is closed under convex combinations
and not the set of ROSt's satisfying GG, we conjecture the following:
\begin{conj}\label{conj:SS->GG}
The laws of the ROSt's satisfying the Ghirlanda-Guerra identities correspond to the extremes of the convex set of laws of the stochastically stable ROSt's.
\end{conj}


\subsection{Application to Spin Glasses}
\label{section: spin glass}
Consider a Gaussian spin glass Hamiltonian $H_N$ and its Gibbs measures $\G_{\beta,N}$ at $\beta$ as defined
in the introduction.  
It is straightforward to construct a ROSt $Q^{\beta,N}$ from $\G_{\beta,N}$ by taking
$$
Q^{\beta,N}=\{ \ r(\sigma^i,\sigma^j)\ \}_{i,j\in\N}
$$
where $(\sigma^i)_{i\in\N}$ are elements of $\{-1,+1\}^N$ sampled iid from $\G_{\beta,N}$. 
Here the form ``$r$'' is the one of the given spin glass model, but the same construction holds for any positive definite symmetric form on $\{-1,+1\}^N$.
Note that the randomness of $Q^{\beta,N}$ comes from the sampling and the randomness of $\G_{\beta,N}$ itself.

The considerations of the previous sections ensure that the random matrix $Q^{\beta,N}$ can be constructed from a random sampling measure $\mu_{\beta,N}$ 
on $\B$. In particular, for any function $F$ on $s$ replicas there is the identity
\begin{equation}
\label{eqn: mu identity}
E \ \G_{\beta,N}^{\otimes s}(F( r(\sigma^k,\sigma^l) ;1\leq k<l\leq s))=\E \ \mu_{\beta,N}^{\otimes s}(F( v^k\cdot v^l; 1\leq k<l\leq s))\ .
\end{equation}
The sequence $(Q^{\beta,N})$ 
in the space of ROSt's has limit points, since the space is compact. 
These limits should retain in some ways the Gibbsian nature of $\G_{\beta,N}$. 
Two properties of the limit points seem to be of importance:
the Ghirlanda-Guerra identities as presented earlier and the stochastic stability for the sequence of sampling measure $(\mu_{\beta,N})$.
Stochastic stability for the sequence defined below does not directly translate into stochastic stability as presented earlier 
of the limit ROSt. As we shall see, Theorem \ref{thm:continuity} fills this gap.

Stochastic stability was first introduced in \cite{AC} and originally defined as follows.
Let $(\G_{\beta,N})$ be a sequence of Gibbs measures at inverse temperature $\beta$.
Consider similarly as before the Gaussian field $(l(\sigma): \sigma\in\{-1,+1\}^N)$ independent of $H_N$, whose law we again denote by $E_l$, with covariance 
$E_l[l(\sigma)l(\sigma')]=r(\sigma,\sigma')$.
For some $\lambda>0$, we consider the mapping \eqref{eqn:mapping}. Since here $\|\sigma\|=1$ for any vector in the support, it reduces to
\begin{equation}
\label{eqn:mapping N}
\G_{\beta,N}(\sigma)\mapsto \frac{\G_{\beta,N}(\sigma)e^{\lambda l(\sigma)}}{\G_{\beta,N}(e^{\lambda l(\sigma)})}\ .
\end{equation}
It is readily checked by property of Gaussians that the image in \eqref{eqn:mapping N} has the same law as the Gibbs measure at temperature $\sqrt{\beta^2+\lambda^2/N}$.
The original idea of \cite{AC} is that a continuous dependence on $\beta$ of the Gibbs measure is equivalent in the limit $N\to\infty$ to stability of the measure under the mapping \eqref{eqn:mapping N}.
\begin{df}
\label{df: SS sequence}
A sequence of Gibbs measures $(\mathcal{G}_N)$ is said to be stochastically stable if for any $\lambda>0$, $s\in\N$ 
and for any continuous function $F$ on $s$ replicas,
$$
\lim_{N\to\infty}\left|E E_l\left[\frac{\mathcal{G}_N^{\otimes s}(F(\ssigma)e^{\lambda l(\sigma^1)} ... e^{\lambda l(\sigma^s)})}{\mathcal{G}_N^{\otimes s}(e^{\lambda l(\sigma^1)} ... e^{\lambda l(\sigma^s)})}\right]
-E \mathcal{G}_N^{\otimes s}(F(\ssigma))\right|=0\ ,
$$
where $F(\ssigma)$ stands for $F(r(\sigma^l,\sigma^k); 1\leq k<l\leq s)$.
\end{df}
It was shown in \cite{CG} that stochastic stability in the sense of Definition \ref{df: SS sequence} holds on average over $\beta$. An improvement was made in \cite{tal_pure} where
it is proved that for any $\beta$ and sequence of Gibbs measures $(\G_{\beta,N})$ there exists a sequence $\beta_N\to\beta$ such that the sequence $(\mathcal{G}_{\beta_N,N})$ is stochastically stable. We show here that stochastic stability holds at any $\beta$ where $\lim_{N\to\infty}\frac{1}{N}E \log Z_N(\beta)$ is differentiable. 
This is an analogue of an elegant result of Panchenko \cite{panchenko_GG2} which proves the validity of the Ghirlanda-Guerra identities under the same hypothesis.

\begin{prop}
\label{thm: diff->SS}
If $\lim_{N\to\infty}\frac{1}{N}E \log Z_N(\beta)$ 
exists and is differentiable at $\beta>0$,
then the sequence of Gibbs measures $(\G_{\beta,N})$ is stochastically stable in the sense of Definition \ref{df: SS sequence}.
\end{prop}
\begin{proof}
Let $F$ be a continuous function on $s$ replicas. We must have $|F|\leq C$ for some $C>0$.
Straightforward differentiation yields
$$
\partial_\beta \ E \G_{\beta,N}^{\otimes s}(F(\ssigma))
= s \ E \G_{\beta,N}^{\otimes s}\Big[\Big(H_N(\sigma^1)-\G_{\beta,N}\big(H_N(\sigma)\big)\big)F(\ssigma)\Big]\ .
$$ 
Hence the absolute value of the derivative is bounded by 
$$s C \ E\G_{\beta,N}\left|H_N(\sigma^1)-\G_{\beta,N}\big(H_N(\sigma)\big)\right|\ .$$
Define $\beta(\lambda):=\sqrt{\beta^2+\lambda^2/N}$. By integration, we have the bound
\begin{multline}
\label{eqn:diff->SS1}
\left|E \mathcal{G}_{\beta(\lambda),N}^{\otimes s }\big(F(\ssigma)\big)- E \G_{\beta,N}^{\otimes s }\big(F(\ssigma)\big)\right|\\
\leq sC \int_\beta^{\beta(\lambda)}E\mathcal{G}_{\beta',N}\left|H_N(\sigma)-\mathcal{G}_{\beta',N}(H_N(\sigma))\right|d\beta'\ .
\end{multline}
It remains to show that the integral goes to zero under the assumption. We have by Cauchy-Schwarz inequality
\begin{multline*}
\int_\beta^{\beta(\lambda)}E\mathcal{G}_{\beta',N}\left|H_N(\sigma)-\mathcal{G}_{\beta',N}(H_N(\sigma))\right|d\beta'
\\\leq
N^{1/2}\big(\beta(\lambda)-\beta\big)^{1/2}\left(\int_\beta^{\beta(\lambda)}\frac{1}{N}E\mathcal{G}_{\beta',N}\left|H_N(\sigma)-\mathcal{G}_{\beta',N}(H_N(\sigma))\right|^2 d\beta'\right)^{1/2}\ .
\end{multline*}
Since $\beta(\lambda)-\beta=\frac{\lambda^2}{2\beta N}+ O(\frac{1}{N^2})$, 
stochastic stability in the sense of Definition \ref{df: SS sequence} 
would follow by equation \eqref{eqn:diff->SS1} if the term in the parentheses divided by $N$ goes to zero as $N\to\infty$. 
Define $f_N(\beta):= \frac{1}{N}E \log Z_N(\beta)$. $f_N$ is famously convex.
It is easily checked by differentiation that the integral to bound is exactly
$$
f_N'(\beta(\lambda))-f_N'(\beta).
$$  
Since $\beta(\lambda)\to\beta$, the above goes to zero by a simple result of convexity (see Lemma \ref{lem:convexity}) as well as the assumptions on the existence of the limit and the derivative at $\beta$.
\end{proof}
A direct consequence of Theorem \ref{thm:continuity} is that stochastic stability as a property of the sequence becomes a property of the limit points.
\begin{cor}
\label{cor: diff->SS}
If $\lim_{N\to\infty}\frac{1}{N}E \log Z_N(\beta)$ 
exists and is differentiable at $\beta>0$,
then the limit points in the space of ROSt's of the sequence $(Q^{\beta,N})$ are stochastically stable in the sense of Definition \ref{df: SS}.
\end{cor}

Another notable application of the continuity of the cavity mapping is
a proof that the Parisi functionals, entering in the Parisi formula of the SK model, are continuous as functionals on ROSt's. 
This is proved in Section \ref{section: parisi functionals} and can be seen as a generalization of a theorem of Guerra that proved
the continuity within the subset of Ruelle Probability Cascades \cite{guerra_bound,arguin_rpc}.
Assuming that the Ultrametricity Conjecture holds, 
this fact, together with Guerra's bound and Corollary \ref{cor: diff->SS},
provides a proof of the Parisi formula for the SK model (and more generally, for
mixed $p$-spin models with even $p$'s). This is done in Section \ref{section: parisi formula}.


\section{The Cavity Mapping}
\label{section:cavity}
\subsection{Continuity}
In this section, we prove the continuity of the mapping \eqref{eqn:mapping}. 
Let $\psi$ be a function in $\mathcal{C}^1(\R)$ 
with bounded derivative, $|\psi'|\leq C$ for some $C>0$. Assume without loss of generality that $\psi(0)=0$.
Let $z$ a standard Gaussian variable with law $P_z$ and expectation $E_z$.
Note that, under these assumptions, $|\psi(z)|\leq C |z|$, and in particular $E_z[e^{m\psi(z)}]<\infty$ for any $m\in\R$. Define
\begin{equation}
\begin{aligned}
\label{eqn:G1}
G(v,z,l):= e^{\psi( l(v)+ z \sqrt{1-\|v\|^2})} \ .
\end{aligned}
\end{equation}
In the cases $\psi(x)=x$ and $\psi(x)=\log \cosh x$, we have the simplification 
$$E_zG(v,z,l)=e^{\psi( l(v))}e^{\frac{1}{2}(1-\|v\|^2)}\ .$$
When dealing with $s$ replicas, that is $\vec{v}=(v^1,...,v^s)\in\B^s$ and $\vec{z}=(z^1,...,z^s)$ where $\vec{z}$ has law $P_z^{\otimes s}$, write
\begin{equation}
\begin{aligned}
\label{eqn:G2}
G(\vec{v},\vec{z},l):=\prod_{i=1}^sG(v^i,z^i,l)\ .
\end{aligned}
\end{equation}
Consider the generalized cavity mapping \eqref{eqn:mapping2}:
\begin{equation}
\mu \mapsto \Phi_{\psi,l}\mu:=\frac{ \mu(dv) \ E_zG(v,z,l) }{\mu(E_zG(v,z,l))}\ ,
\end{equation}
where $l(v)$ is the cavity field with covariance $v\cdot v'$. The proof of the continuity is the same
when the covariance of the field is $c(v\cdot v')$, for a suitable function $c$, if $\|v\|^2$ in the definition of $G$
is replaced by $c(\|v\|^{2})$. We restrict ourselves to the case where $c$ is the identity to simplify notation.

Consider a function $F=F(\vec{v},\vec{z},l)$ that depends measurably on $s$ replicas,
the Gaussian vector $\vec{z}$ and the cavity field. The main ingredient of the proof of continuity is an approximation
in the spirit of the weak law of large numbers.
\begin{lem}
\label{lem:weak}
Let $\mu\in\M(\B)$. Let $F(\vec{v},\vec{z},l)$ be 
a function such that $\vec{v}\mapsto E_lE_z^{\otimes s}[F(\vec{v},\vec{z},l)^2]$ is a continuous function on $s$ replicas. Let $(\vec{v}^r,\vec{z}^r)$, $r=1,...,n$, be  $n$ independent copies of $(\vec{v},\vec{z})$ 
sampled from $(\mu\times P_z)^{\otimes s}$. Then
$$
E_l \ (\mu\times E_z)^{\otimes ns}\left(\frac{1}{n}\sum_{r=1}^n F(\vec{v}^{r},\vec{z}^r,l)
-(\mu\times E_z)^{\otimes s}\Big(F(\vec{v},\vec{z},l)\Big)\right)^2
\leq \frac{C}{n}
$$
where $C$ is a positive constant that depends on $F$ but not on $\mu$.
\end{lem}
\begin{proof}
Using the fact that under $(\mu\times E_z)^{\otimes ns}$, the variables $F(\vec{v}^{r},\vec{z}^r,l)$, $r=1,...,n$ are i.i.d., we have
$$
\begin{aligned}
& (\mu\times E_z)^{\otimes ns}\left(\frac{1}{n}\sum_{r=1}^n F(\vec{v}^{r},\vec{z}^r,l)
-(\mu\times E_z)^{\otimes s}\Big(F(\vec{v},\vec{z},l)\Big)\right)^2\\
&=\frac{1}{n} (\mu\times E_z)^{\otimes s}\Big(F(\vec{v},\vec{z},l)^2\Big)-
\frac{1}{n}\Big( (\mu\times E_z)^{\otimes s}\big(F(\vec{v},\vec{z})\big)\Big)^2\\
&\leq \frac{1}{n} (\mu\times E_z)^{\otimes s}\Big(F(\vec{v},\vec{z},l)^2\Big)\ .
\end{aligned}
$$
Now integrate over the $l$'s on both sides and use Fubini's theorem on the right-hand-side.  
By assumption, $\vec{v}\mapsto E_lE_z^{\otimes s}[F(\vec{v},\vec{z},l)^2]$ is a bounded function on $\B^s$. The conclusion follows from this.
\end{proof}
It will be useful to impose stronger conditions on $F(\vec{v},\vec{z},l)$, namely that
$F(\vec{v},\vec{z},l)$ is a function of the variables $ l(v^i)+z^i\sqrt{1-\|v^i\|^2}$, $i=1,...,s$:
\begin{equation}
\label{eqn:ass F}
\begin{aligned}
F(\vec{v},\vec{z},l)=F(l(v^i)+z^i\sqrt{1-\|v^i\|^2}, i=1,...,s)\ ;
\end{aligned}
\end{equation}
and that for any $m\in\R$,
\begin{equation}
\label{eqn:ass2 F}
\begin{aligned}
\vec{v}\mapsto E_lE_z^{\otimes s}[F(\vec{v},\vec{z},l)^m]\text{ is a continuous function on $s$ replicas.}
\end{aligned}
\end{equation}

\begin{lem}
\label{lem:bound}
The function $G(\vec{v},\vec{z},l)$ satisfies \eqref{eqn:ass F} and \eqref{eqn:ass2 F}.
\end{lem}
\begin{proof}
The dependence in $G$ is uniquely on the Gaussian variables $l(v^i)+z^i\sqrt{1-\|v\|^2}$, $i=1,...,s$.
These variables have variance $1$ and covariance $v^i\cdot v^j$, $i\neq j$. Therefore the function 
$\vec{v}=(v^1,...,v^s)\mapsto  E_lE_z^{\otimes s}[G(\vec{v},\vec{z},l)^m]$ depends
only on the inner products between distinct replicas $v^1$,...,$v^s$. By writing down the Gaussian integral
$E_lE_z^{\otimes s}$ explicitly,
we see that to prove it is continuous on $\B^s$ it suffices 
to show that, for any $m\in\R$, the function
$\vec{v} \mapsto E_lE_z^{\otimes s}[G(\vec{v},\vec{z},l)^m]$ is bounded uniformly in $\vec{v}$.
By H\"older's inequality,
$$
E_lE_z^{\otimes s}\left[\prod_{i=1}^sG(v^i,z^i,l)^{m}\right]\leq \prod_{i=1}^s \left\{E_lE_{z^i}[G(v^i,z^i,l)^{ms}]\right\}^{1/s}=E_z[e^{sm\psi(z)}]<\infty, 
$$
which shows the desired uniform bound.
\end{proof}

Theorem \ref{thm:continuity} is a straightforward consequence of the following result with $F(\vec{v},\vec{z},l)=F(\vec{v})$.
The idea of the proof is to linearize functions involving the measure $\mu\times E_z$ using the empirical approximation of Lemma 
\ref{lem:weak}.
\begin{thm}\label{thm2:continuity}
Let $F(\vec{v},\vec{z},l)$ be a function satisfying \eqref{eqn:ass F} and \eqref{eqn:ass2 F}  and let $G(\vec{v}, \vec{z},l)$ be as in \eqref{eqn:G2}.
Then the mapping
\begin{equation}
\label{eqn:mapping G}
\prob \mapsto 
\E E_l\left[ \frac{(\mu\times E_z)^{\otimes s}\Big(F(\vec{v},\vec{z},l)G(\vec{v},\vec{z},l)\Big)}{(\mu\times E_z)^{\otimes s}\Big(G(\vec{v},\vec{z},l)\Big)}\right]
\end{equation}
is a continuous functional on the space of ROSt's.
\end{thm}

\begin{proof}
We rearrange the image of the mapping,
\begin{equation}
\label{eqn:proof}
\begin{aligned}
&\E E_l\left[ \frac{(\mu\times E_z)^{\otimes s}\Big(F(\vec{v},\vec{z},l)G(\vec{v},\vec{z},l)\Big)}{(\mu\times E_z)^{\otimes s}\Big(G(\vec{v},\vec{z},l)\Big)}\right]=
\\
&\hspace{2cm} \E E_l\ (\mu\times E_z)^{\otimes ns}\left[\frac{\frac{1}{n}\sum_{r=1}^n F(\vec{v}^{r},\vec{z}^r,l)G(\vec{v}^r,\vec{z}^r,l)+
d_\mu(FG)}{\frac{1}{n}\sum_{r=1}^n G(\vec{v}^r,\vec{z}^r,l)+d_\mu(G)}\right]
\end{aligned}
\end{equation}
where
$$
\begin{aligned}
d_\mu(FG)&:=(\mu\times E_z)^{\otimes s}\Big(F(\vec{v},\vec{z},l)G(\vec{v},\vec{z},l)\Big)-\frac{1}{n}\sum_{r=1}^n F(\vec{v}^{r},\vec{z}^r,l)G(\vec{v}^{r},\vec{z}^r,l)\\
d_\mu(G)&:=(\mu\times E_z)^{\otimes s}\Big(G(\vec{v},\vec{z},l)\Big)-\frac{1}{n}\sum_{r=1}^n G(\vec{v}^{r},\vec{z}^r,l)\ .
\end{aligned}
$$
We introduce the empirical average of $f$ as a function on the $n$ copies $(\vec{v}^r,\vec{z}^r)_{r=1}^n$ and the cavity field
\begin{equation}
\label{eqn:F tilde}
\tilde{F}((\vec{v}^r,\vec{z}^r)_r;l):=\frac{\frac{1}{n}\sum_{r=1}^n F(\vec{v}^{r},\vec{z}^r,l)G(\vec{v}^{r},\vec{z}^r,l)}{\frac{1}{n}\sum_{r=1}^n G(\vec{v}^{r},\vec{z}^r,l)}\ .
\end{equation}
In this notation, elementary manipulations of \eqref{eqn:proof} give
\begin{equation}
\label{eqn: cont eqn}
\begin{aligned}
&\E E_l\left[ \frac{(\mu\times E_z)^{\otimes s}\Big(F(\vec{v},\vec{z},l)G(\vec{v},\vec{z},l)\Big)}{(\mu\times E_z)^{\otimes s}\Big(G(\vec{v},\vec{z},l)\Big)}\right]
-\E E_l\ (\mu\times E_z)^{\otimes ns}\left(\tilde{F}((\vec{v}^r,\vec{z}^r)_r;l)\right)
\\ &\hspace {3cm}=
 \E E_l\left[\frac{(\mu\times E_z)^{\otimes ns}\Big(d_\mu(FG)-\tilde{F}((\vec{v}^r,\vec{z}^r)_r;l)\ d_\mu(G)\Big)}
{(\mu\times E_z)^{\otimes s}\Big(G(\vec{v},\vec{z},l)\Big)}\right]\ .
\end{aligned}
\end{equation}
We now bound the absolute value of the difference appearing in the left-hand side of \eqref{eqn: cont eqn} uniformly in the possible laws $\prob$. 
We write for simplicity $\gamma=\sqrt{2}$ and $\bar{\gamma}=(1-\frac{1}{\gamma})^{-1}=\frac{\sqrt{2}}{\sqrt{2}-1}$. 
We use H\"older's inequality followed by Jensen's inequality and Fubini's theorem to get the following upper bound of the right-hand side:
\begin{multline*}
\left\{\E\mu^{\otimes ns}\Big( E_lE_z^{\otimes ns}|d_\mu(FG)-\tilde{F}((\vec{v}^r,\vec{z}^r)_r;l)\ d_\mu(G)|^\gamma\Big)\right\}^{1/\gamma}\\
\times 
\left\{\E\mu^{\otimes ns}\Big( E_lE_z^{\otimes ns}[G(\vec{v},\vec{z},l)^{-\bar{\gamma}}]\Big)\right\}^{1/\bar{\gamma}}\ .
\end{multline*}
The second term is bounded by a constant uniform in the laws $\prob$ by Lemma~\ref{lem:bound}.
(Essentially, the term with negative power $-\bar{\gamma}$ is bounded because the Laplace functional of a Gaussian variable is
well-defined everywhere.)
As for the first term, by the triangle inequality, it is smaller than
$$
\begin{aligned}
&\left\{\E\mu^{\otimes ns}\Big( E_lE_z^{\otimes ns} |d_\mu(FG)|^{\gamma}\Big)\right\}^{1/\gamma}\\
&\hspace{4cm}+
\left\{\E\mu^{\otimes ns}\Big( E_lE_z^{\otimes ns}|\tilde{F}((\vec{v}^r,\vec{z}^r)_r;l)\ d_\mu(G)|^{\gamma}\Big)\right\}^{1/\gamma}.
\end{aligned}
$$
The first term is bounded by the same expression with $\gamma$ replaced by $2$. 
By Lemma \ref{lem:weak}, it is thus bounded by $C/\sqrt{n}$ uniformly in $\prob$.
The second term is bounded similarly after an application of H\"older's inequality with $\gamma$.
We note that, in this case, the resulting term 
$\E\mu^{\otimes ns}\Big( E_lE_z^{\otimes ns}|\tilde{F}((\vec{v}^r,\vec{z}^r)_r;l)|^{\gamma\bar{\gamma}}\Big)$
is bounded uniformly in the $\prob$'s since
\begin{equation}
\label{eqn: estimate F}
\begin{aligned}
E_lE_z^{\otimes ns}[\tilde{F}((\vec{v}^r,\vec{z}^r)_r;l)^{\gamma\bar{\gamma}}]
\leq &\Big\{E_l E_z^{\otimes s}\left[F(\vec{v},\vec{z},l)^{3\gamma\bar{\gamma}}\right]\Big\}^{1/3}  \left\{E_l E_z^{\otimes s}\left[G(\vec{v},\vec{z},l)^{3\gamma\bar{\gamma}}\right]\right\}^{1/3} \\
&\qquad \qquad \times 
\Big\{E_lE_z^{\otimes s}\left[ G(\vec{v},\vec{z},l)^{-3\gamma\bar{\gamma}}\right]\Big\}^{1/3}\ ,
\end{aligned}
\end{equation}
where we have used H\"older's inequality and Jensen's inequality.
The right-hand side is bounded because of \eqref{eqn:ass2 F} and Lemma  \ref{lem:bound}.
Putting all this together, we can write from \eqref{eqn: cont eqn} for a possibly different $C$, not depending on $\prob$, 
\begin{align}
\biggl|\E E_l\left[ \frac{(\mu\times E_z)^{\otimes s}\Big(F(\vec{v},\vec{z},l)G(\vec{v},\vec{z},l)\Big)}{(\mu\times E_z)^{\otimes s}\Big(G(\vec{v},\vec{z},l)\Big)}\right]
&-  \E \mu^{\otimes ns}\Big(E_lE_z^{\otimes ns}[\tilde{F}((\vec{v}^r,\vec{z}^r)_r;l)]\Big)\biggr| \nonumber\\
&\leq  \frac{C}{\sqrt{n}}.
\label{eqn:approx}
\end{align}
The above approximation is useful because it is uniform in the laws $\prob$, and it linearizes the dependence on $\mu$
(allowing the use of Fubini's theorem to exchange the integration $E_l$ and $\mu$).

If $\E$ and $\E'$ are two ROSt's with sampling measures $\mu$ and $\mu'$, we have by \eqref{eqn:approx}
$$
\begin{aligned}
&\left|\E E_l\left[ \frac{(\mu\times E_z)^{\otimes s}\Big(F(\vec{v},\vec{z},l)G(\vec{v},\vec{z},l)\Big)}{(\mu\times E_z)^{\otimes s}\Big(G(\vec{v},\vec{z},l)\Big)}\right]\right.
\\
& \qquad \qquad \left.- \E' E_l\left[ \frac{(\mu'\times E_z)^{\otimes s}\Big(F(\vec{v},\vec{z},l)G(\vec{v},\vec{z},l)\Big)}{(\mu'\times E_z)^{\otimes s}\Big(G(\vec{v},\vec{z},l)\Big)}\right]
\right|
\\
&\leq 
\left|\E\mu^{\otimes ns}\Big(E_lE_z^{\otimes ns}[\tilde{F}((\vec{v}^r,\vec{z}^r)_r;l)]\Big)-
\E'\mu'^{\otimes ns}\Big(E_lE_z^{\otimes ns}[\tilde{F}((\vec{v}^r,\vec{z}^r)_r;l)]\Big)\right| \\
&\qquad \qquad  + \frac{2C}{\sqrt{n}}.
\end{aligned}
$$
Since $C$ is uniform in the choice of $\prob$, the continuity will follow if
$$(\vec{v}^r)_{r=1}^n\mapsto E_lE_z^{\otimes ns}[\tilde{F}((\vec{v}^r,\vec{z}^r)_r;l)] $$
is a continuous function on $ns$ replicas.
We have that
$$
E_lE_z^{\otimes ns}[\tilde{F}((\vec{v}^r,\vec{z}^r)_r;l)]=
E_lE_z^{\otimes ns}\left[\frac{\frac{1}{n}\sum_{r=1}^nF(\vec{v}^r,\vec{z}^r,l)G(\vec{v}^{r},\vec{z}^r,l)}{\frac{1}{n}\sum_{r=1}^n G(\vec{v}^{r},\vec{z}^r,l)}\right].
$$
But the integration is on the $ns$ Gaussian variables of the form $l(v)+z\sqrt{1-\|v\|^2}$ that have variance $1$
and covariance $v\cdot v'$.
Therefore the expectation is a function on $ns$ replicas, since it depends only on the inner products between distinct replicas and
not on their norm. It is readily seen from writing the Gaussian integral that the continuity will follow
if the expectation is bounded as a function on $\B^{ns}$. But this is a consequence of \eqref{eqn: estimate F} taking $\gamma\bar\gamma=1$.
\end{proof}

Compositions of cavity mappings take a simple form when $\psi$ is linear.
This will be needed in Section \ref{section:properties}.
\begin{lem}
\label{lem: composition}
Let $\lambda,\lambda'>0$, $c$ be as in \eqref{eqn:c} and $\prob$ be the law of a ROSt. If $\prob'$ is the image of $\prob$ under $\Phi_{\lambda,l_c}$
and $\prob''$ is the image of $\prob'$ under $\Phi_{\lambda',l'_c}$, then $\prob''$ is the image of $\prob$ under $\Phi_{\sqrt{\lambda^2+\lambda'^2},l_c}$.
In particular, if $\prob$ is invariant under $\Phi_{\lambda,l_c}$, then it is invariant under $\Phi_{\sqrt{T}\lambda,l_c}$ for any $T\in\N$.
\end{lem}
\begin{proof}
The second statement follows directly from the first. 
Let $\mu$ be the sampling measure of $\prob$.
Then $\mu'=\Phi_{\lambda,l_c}\mu$ is the sampling measure of $\prob'$.
By the definition \eqref{eqn:mapping2} of the mapping,
$$
\Phi_{\lambda',l_c'}\mu'=\frac{e^{\lambda' l_c'( v)-\frac{\lambda'^2}{2}\| v\|^2}\ \mu'(d v) }{\mu'(e^{\lambda' l'_c( v)-\frac{\lambda'^2}{2}\| v\|^2})}\ ,
$$
and by the definition of $\mu'$,
$$
\Phi_{\lambda',l_c'}\mu'=\frac{e^{\lambda' l_c'( v)+\lambda l_c( v)-\frac{\lambda'^2+\lambda^2}{2}\| v\|^2}\ \mu(d v) }{\mu(e^{\lambda' l_c'( v)+\lambda l_c( v)-\frac{\lambda'^2+\lambda^2}{2}\| v\|^2})}\ .
$$
The claim follows from the fact that $\lambda'l_c'+\lambda l_c \distrib \sqrt{\lambda'^2+\lambda^2} ~ l_c$.
\end{proof}


\subsection{The Parisi Functionals}
\label{section: parisi functionals}
Theorem \ref{thm2:continuity} that proves the continuity of the cavity mapping also provides continuity for an important class of functionals on ROSt's, 
the so-called {\it Parisi functionals}.
They enter in particular in the Parisi formula for the free energy of the SK model.

Let $\psi\in C^1(\R)$ with bounded derivative and with $\psi(0)=0$.
Consider a cavity field $l_c$ with covariance 
$E_l[l_c(v)l_c(v')]=c(v\cdot v')$ for some function $c$ as in \eqref{eqn:c}.
Define the Parisi functional for a ROSt $\prob$ and parameters $\lambda\geq 0$ as
\begin{eqnarray}
\label{eqn:Parisi functional}
\mathcal{P}(\lambda,\prob):=\E \ E_l \left[\log \ (\mu\times E_z) \Big(e^{\psi(\lambda l_c(v)+\lambda z\sqrt{1-c(\|v\|^2)})}\Big)\right]
\end{eqnarray} 
where $z$ is a standard Gaussian variable with expectation $E_z$.
We drop the dependence on $c$ and on $\psi$ in the notation for simplicity. 
The cases where $\psi(x)=\log\cosh x$ or $\psi(x)=x$,
and
$$
\begin{aligned}
\lambda^2=\sum_{p\geq 0}\lambda_p^2 \qquad c(x)= \frac{1}{\lambda^2}\sum_{p\geq 0}\lambda_p^2 ~ x^{p}
\end{aligned}
$$
will play an important role in Section \ref{section: parisi formula}. In these cases, the variable $z$ can 
be integrated. The functional \eqref{eqn:Parisi functional} reduces to
\begin{equation}
\label{eqn:Parisi functional p}
\mathcal{P}(\vec{\lambda},\prob):=\frac{\lambda^2}{2} +\E \ E_l \left[\log \mu \Big(e^{\psi(\sum_{p\geq 0}\lambda_p l_p(v))-\frac{1}{2}\sum_{p\geq 0}\lambda_p^2\|v\|^{2p}}\Big)\right]\ ,
\end{equation}
where $\vec{\lambda}=(\lambda_p)_{p\geq 0}$ and $E_l$ is the expectation over independent cavity fields $l_p$, $p\in\N\cup\{0\}$.
The dependence on each $\lambda_p$ is emphasized here by writing $\mathcal{P}(\vec{\lambda},\prob)$ instead of $\mathcal{P}(\lambda,\prob)$.

Continuity of the Parisi functional on the subset of Ruelle Probability Cascades has
been established by Guerra \cite{guerra_bound}. In this case, the Parisi functionals equals $f(0,0)$, where $f(q,y)$ is the solution of a p.d.e. with
final condition $f(1,y)=\psi(\lambda y)$, see e.g. \cite{ASS,arguin_rpc} for details. 
It turns out that continuity properties of the Parisi functionals that are needed in the proof of the Parisi formula can be established for the whole space of ROSt's.

\begin{prop}
\label{prop:continuity parisi}
Let $\mathcal{P}(\lambda, \prob)$ be as in \eqref{eqn:Parisi functional}. 
Then, for any $\psi\in C^1(\R)$ with bounded derivative and $c$ as in \eqref{eqn:c},
\begin{enumerate}
\item $\prob\mapsto \mathcal{P}(\lambda,\prob)$ is a continuous functional on ROSt's for every $\lambda\geq0$.
\item In the cases where $\psi(x)=\log \cosh x$ and $\psi(x)=x$ as in \eqref{eqn:Parisi functional p}, it holds for any $\delta>0$ and any given $\vec{\lambda'}=(\lambda'_p)_{p\geq 0}\in l^2$ that
$$
\text{ if $\|\vec{\lambda} -\vec{\lambda}'\|_{2}<\delta$, then }\max_{\prob \text{ ROSt}} \left|\mathcal{P}(\vec{\lambda}, \prob) -  \mathcal{P}(\vec{\lambda}',\prob)\right|\leq K \|\vec{\lambda}-\vec{\lambda}'\|_{2}
$$
where $K>0$ depends on $\vec{\lambda}'$ and $\delta$.
\end{enumerate}
\end{prop}
\begin{proof}
We prove the first claim. 
We show that $\mathcal{P}(\lambda, \prob)$ is differentiable in $\lambda$ and 
that the derivative is a continuous function on ROSt's. The conclusion follows by integration.

We write for short
$$
g(v,z):=l(v)+z\sqrt{1-c(\|v\|^2)}\ ,
$$
omitting the dependence on $l$. Note that under $E_lE_z$, the $g$'s are all centered Gaussian variables
of variance $1$.
Suppose first that the derivative can be taken inside the expectations $\E \ E_l$ and $\mu\times E_z$, then
one would get
\begin{equation}
\label{eqn:deriv}
\partial_\lambda \mathcal{P}(\lambda, \prob)= 
\E E_l\left[\frac{(\mu\times E_z)
\Big(g(v,z)\psi'(\lambda g(v,z))e^{\psi(\lambda g(v,z))}
\Big)}
{(\mu\times E_z)\Big(e^{\psi(\lambda g(v,z))}\Big)}\right]\ .
\end{equation}
By Theorem \ref{thm2:continuity}, the derivative will be continuous as a functional on ROSt's if $g(v,z)\psi'(\lambda g(v,z))$ satisfies \eqref{eqn:ass F}
and \eqref{eqn:ass2 F}.
This is straightforward since $\psi'$ is bounded. 

It remains to show \eqref{eqn:deriv}. The derivative can be taken inside the expectation $\E E_l$ if we show that for $\delta$ small enough
\begin{equation}
\label{eqn:unif int log}
\sup_{|\lambda-\lambda'|<\delta} \frac{1}{(\lambda-\lambda')^2}\ 
\E E_l \left(\log 
\frac{(\mu\times E_z)\Big(e^{\psi(\lambda g(v,z))}\Big)}
{(\mu\times E_z)\Big(e^{\psi(\lambda' g(v,z))}\Big)}\right)^2<\infty\ .
\end{equation}
This effectively demonstrates the uniform integrability of the collection of random variables indexed by $\lambda$ with $|\lambda-\lambda'|< \delta$ and given by
$$\frac{1}{\lambda-\lambda'}\Big|\log 
(\mu\times E_z)\Big(e^{\psi(\lambda g(v,z))}\Big)-\log
(\mu\times E_z)\Big(e^{\psi(\lambda' g(v,z))}\Big)\Big|\ .
$$
Write $R(\lambda,\lambda')$ for the random variable inside the $\log$ in \eqref{eqn:unif int log}. We split the integration for the events $\{R(\lambda,\lambda')<1\}$ and $\{R(\lambda,\lambda')>1\}$.
The expectation in \eqref{eqn:unif int log} becomes
\begin{align*}
& \E E_l \left[\left(\log R(\lambda,\lambda')\right)^2\ ; R(\lambda,\lambda')>1\right]
\\
&\qquad \qquad + \E E_l \left[\left(\log R(\lambda,\lambda')^{-1}\right)^2;R(\lambda,\lambda')^{-1}>1\right].
\end{align*}
Since $\log x\leq x-1$ for $x>1$, we can bound the above by
\begin{multline}
 \E E_l \left(\frac{(\mu\times E_z)\Big(e^{\psi(\lambda g(v,z))}-e^{\psi(\lambda' g(v,z))}\Big)}
 {(\mu\times E_z)\Big(e^{\psi(\lambda g(v,z))}\Big)}\right)^2+\\
  \E E_l \left(\frac{(\mu\times E_z)\Big(e^{\psi(\lambda g(v,z))}-e^{\psi(\lambda' g(v,z))}\Big)}
 {(\mu\times E_z)\Big(e^{\psi(\lambda' g(v,z))}\Big)}\right)^2
\ .
\end{multline}
Using successively Cauchy-Schwarz inequality followed by Jensen's inequality on each term yields the upper bound
\begin{multline}
\label{eqn:unif bound}
\left\{ \E E_l \left[(\mu\times E_z)\Big( e^{\psi(\lambda g(v,z))}-e^{\psi(\lambda' g(v,z))} \Big)^4\right]\right\}^{1/2} \\
\times \left\{\left( \E E_l\left[(\mu\times E_z)\Big(e^{-4\psi(\lambda g(v,z))}\Big)\right]\right)^{1/2}\right.\\
\qquad +
 \left. \left(\E E_l\left[(\mu\times E_z)\Big(e^{-4\psi(\lambda' g(v,z))}\Big)\right]\right)^{1/2}\right\}
\end{multline}
The terms in the second bracket can be integrated by Fubini's and since the fields $g(v,z)$ all have variance $1$ under $E_lE_z$, 
the dependence on $\E$ drops leaving
$$
E_z\left[e^{-4\psi(\lambda z)}\right]^{1/2}+ E_z\left[e^{-4\psi(\lambda' z)}\right]^{1/2}\ .
$$
The fact that $|\psi'|\leq C$ and $\psi(0)=0$ implies that $|\psi(\lambda z)|\leq \lambda C |z|$.
And since $|\lambda-\lambda'|<\delta$,  $|\psi(\lambda z)|\leq C(\lambda'+\delta)|z|$.
Hence after integration, the term in the parentheses is seen to be uniformly bounded by
a term that depends only on $\lambda'$ and $\delta$.
It remains to bound the first term. After a use of Fubini's theorem, since the fields $g(v,z)$ have variance $1$, it is simply 
\begin{equation}
\label{eqn:unif bound 2}
  \left\{E_z\left(e^{\psi(\lambda z)}-e^{\psi(\lambda' z)}\right)^4\right\}^{1/2}\ .
\end{equation}
Taylor's expansion around $\psi(\lambda' z)$ gives for a certain $\bar{z}$ between $\lambda' z$ and $\lambda z$,
\begin{equation}
\label{eqn:expansion}
e^{\psi(\lambda z)}-e^{\psi(\lambda' z)}=(\lambda-\lambda')z \psi'(\bar{z})e^{\psi(\bar{z})}.
\end{equation}
The latter is bounded above by $C|\lambda-\lambda'| \ |z| e^{C |\bar{z}|}$. The assumption on $\bar{z}$
readily implies that $|\bar{z}|\leq |z| (\lambda' +\delta)$. 
Plugging all this back into \eqref{eqn:unif bound 2}, we get a bound after integration of $z$ which is $(\lambda-\lambda')^2$
times a
constant that depends only $\lambda'$ and $\delta$ . 
Since the term $(\lambda-\lambda')^2$ cancels with the one in \eqref{eqn:unif int log},
it follows that the differential operator can be passed through $\E E_l$. 

To prove \eqref{eqn:deriv}, 
it remains to show
\begin{equation}
\label{eqn:deriv2}
\partial_\lambda \log 
(\mu\times E_z)\Big(e^{\psi(\lambda g(v,z))}\Big)= \frac{(\mu\times E_z)
\Big(g(v,z)\psi'(\lambda g(v,z))e^{\psi(\lambda g(v,z))}
\Big)}
{(\mu\times E_z)\Big(e^{\psi(\lambda g(v,z))}\Big)}\ .
\end{equation}
This is straightforward once it is established that
$$
\partial_\lambda
(\mu\times E_z)\Big(e^{\psi(\lambda g(v,z))}\Big)= (\mu\times E_z)
\Big(g(v,z)\psi'(\lambda g(v,z))e^{\psi(\lambda g(v,z))}
\Big)\ .
$$
To do so, it suffices to justify taking the derivative inside $\mu\times E_z$ by showing
\begin{equation}
\sup_{|\lambda-\lambda'|<\delta} \frac{1}{(\lambda-\lambda')^2}
(\mu\times E_z)\Big(e^{\psi(\lambda g(v,z))}-e^{\psi(\lambda' g(v,z))}\Big)^2<\infty\ .
\end{equation}
Using the expansion \eqref{eqn:expansion} with $g$, the above is bounded by $$C^2 (\mu\times E_z)(g^2 e^{2C(\lambda'+\delta)|g|}).$$
The expectation is finite $P_l$-almost surely as can be seen by integrating over $l$. \eqref{eqn:deriv2} follows. 

To prove the second claim of the proposition, it suffices to show
\begin{equation}
\label{eqn: to bound p}
\E E_l \left| \log\frac{ \mu \Big(e^{\psi(\sum_p\lambda_p l_p(v))-\frac{1}{2}\sum_p\lambda_p^2\|v\|^{2p}}\Big)} {\mu \Big(e^{\psi(\sum_p\lambda'_p l_p(v))-\frac{1}{2}\sum_p\lambda'^2_p\|v\|^{2p}}\Big)}\right|
<K \| \vec{\lambda}-\vec{\lambda}'\|_{2}\ ,
\end{equation}
since the functional $\vec{\lambda}\mapsto \frac{\lambda^2}{2}$ in \eqref{eqn:Parisi functional p}
obviously obeys a similar bound for $\|\vec{\lambda}-\vec{\lambda}'\|_{2}<\delta$.
We write $K$ for a generic term that depends on $\vec{\lambda}'$ and $\delta$ only.
Bounding the logarithm as it was done after \eqref{eqn:unif int log}, one gets that \eqref{eqn: to bound p} is smaller or equal than a term bounded on 
$\{\vec{\lambda}:\|\vec{\lambda}-\vec{\lambda}'\|_{2}<\delta\}$ times 
$$
\left\{ \E E_l ~ \mu \Big(e^{\psi(\sum_p\lambda_p l_p(v))}e^{-\frac{1}{2}\sum_p\lambda_p^2\|v\|^{2p}}-e^{\psi(\sum_p\lambda'_p l_p(v))}e^{-\frac{1}{2}\sum_p\lambda'^2_p\|v\|^{2p}}\Big)^2 \right\}^{1/2}\ .
$$
By adding and subtracting $e^{\psi(\sum_p\lambda'_p l_p(v))}e^{-\frac{1}{2}\sum_p\lambda_p^2\|v\|^{2p}}$ and using the triangle inequality, this is smaller than
$$
\begin{aligned}
\left\{ \E E_l ~ \mu \Big(e^{-\sum_p\lambda_p^2\|v\|^{2p}}\big(e^{\psi(\sum_p\lambda_p l_p(v))}-e^{\psi(\sum_p\lambda'_p l_p(v))}\big)^2\Big) \right\}^{1/2}
+\\
\left\{ \E E_l ~ \mu \Big(e^{2\psi(\sum_p\lambda'_p l_p(v))}\big(e^{-\frac{1}{2}\sum_p\lambda_p^2\|v\|^{2p}}-e^{-\frac{1}{2}\sum_p\lambda'^2_p\|v\|^{2p}}\big)^2\Big) \right\}^{1/2}\ .
\end{aligned}
$$
The terms $e^{-\sum_p\lambda_p^2\|v\|^{2p}}$ and $E_l~ e^{2\psi(\sum_p\lambda'_p l_p(v))}$ are uniformly bounded on $\B$ and on $\{\vec{\lambda}: \|\vec{\lambda}-\vec{\lambda}'\|_{2}<\delta\}$.
Moreover, for $\psi=\log\cosh x$ or $\psi(x)=x$, one has $(e^{\psi(x)}-e^{\psi(y)})^2\leq (x-y)^2e^{2|x|+2|y|}$.
Therefore, applying Cauchy-Schwarz inequality, the claimed bound follows for the first term from the fact that 
$\lambda_p l_p(v)$ has bounded exponential moments on $\B$ and on $\{\vec{\lambda}:\|\vec{\lambda}-\vec{\lambda}'\|_{2}<\delta\}$,
and that 
$$
E_l \big(\sum_p (\lambda_p-\lambda_p')l_p(v)\big)^4=3\left(\sum_p (\lambda_p-\lambda_p')^2\|v\|^{2p}\right)^2\leq 3\|\vec{\lambda}-\vec{\lambda}'\|^4_{2}\ .
$$
As for the second term, the bound follows from 
$$
|e^{-\frac{1}{2}\sum_p\lambda_p^2\|v\|^{2p}}-e^{-\frac{1}{2}\sum_p\lambda'^2_p\|v\|^{2p}}|
\leq \frac{1}{2}\sum_p \lambda_p^2-\lambda_p'^2\leq K \|\vec{\lambda}-\vec{\lambda}'\|_{2}\ .
$$

 \end{proof}


\section{The Parisi Formula from the Ultrametricity Conjecture}
\label{section: parisi formula}

In this section, we show that the Ultrametricity Conjecture \ref{conj:ultra} together with Corollary \ref{cor: diff->SS} and Guerra's upper bound \cite{guerra_bound} yield
an alternative proof of the Parisi formula for the SK model (and more generally, for mixed $p$-spin models with even $p$'s). 
The main new result is Theorem \ref{thm: parisi 1}, which gives a lower bound for the free energy in terms of a functional 
of a limit point of the sequence of Gibbs measures in the space of ROSt's.
The proof uses the treatment of Talagrand for the high-temperature regime (Theorem 2.4.19 in \cite{tal_book})
and the continuity of the Parisi functionals on the space of ROSt's. The compactness
of the set of ROSt's that are invariant under a family of cavity mappings, which is a consequence of the continuity
of such mappings, will also be crucial (see Corollary \ref{cor: beta positive}).
The Parisi formula is then a consequence of
Guerra's bound and the Ultrametricity Conjecture. 

The $p$-spin Hamiltonian $H_N^p$ on $\{-1,+1\}^N$, 
$$(H_{N}^p(\sigma), \sigma\in\{-1,+1\}^N)\ ,$$
is a centered Gaussian vector with covariances
$E H_{N}^p(\sigma)H_{N}^p(\sigma')=N R(\sigma,\sigma')^p$,
where $R(\sigma,\sigma')=\frac{1}{N}\sum_{i=1}^N\sigma_i\sigma_i'$.
Such a process can be represented as
\begin{equation}
\label{eqn: H_p}
H_N^p(\sigma)= \frac{1}{N^{(p-1)/2}}\sum_{i_1,...,i_p=1}^N g_{i_1,...,i_p}\sigma_{i_1}...\sigma_{i_p} \ ,
\end{equation}
where $g_{i_1,...,i_p}$ are standard Gaussians that are independent for distinct indices $(i_1,...,i_p)$.
In this section, we will be concerned with the mixed Hamiltonian
$$
H_{N}(\sigma)=\sum_{p}\beta_p H_N^p(\sigma)
$$
where the summation is restricted throughout this section to $p=1$ and even $p$'s. Take $\beta_p\in\R_+$ for simplicity since $H_N$ has the same law for positive and negative values of the parameters.
For the appropriate cavity field to be well-defined, the parameters $\bbeta:=(\beta_p)$ will need to decay so that
\begin{equation}
\label{eqn: beta}
\sum_p 2^p \ \beta_p^2<\infty \ .
\end{equation}
(We stress that this decay is most likely not optimal.) The $\bbeta$'s satisfying \eqref{eqn: beta} form a Hilbert space under the inner product
$\bbeta\cdot \bbeta'=\sum_p 2^p \beta_p\beta_p'$. The norm is then $\|\bbeta\|^2=\sum_p 2^p \beta_p^2$.

The Gibbs measure at the value $\bbeta$ of the parameters  is
$$
\G_{\bbeta,N}(\sigma)=\frac{\exp{ H_N(\sigma)}}{Z_N(\bbeta)}\ ,
$$
where $Z_N(\bbeta)=\sum_{\sigma}\exp H_N(\sigma)$.
As mentioned in Section \ref{section: spin glass}, 
a ROSt can be constructed from the random probability measures $(\G_{\bbeta,N})$.
The corresponding sequence of laws will be denoted by $(\prob_{\bbeta,N})$. 
The corresponding sampling measures defined through \eqref{eqn: mu identity} will be denoted by $(\mu_{\bbeta,N})$.
Since the space of ROSt's is compact, the sequence $(\prob_{\bbeta,N})$ has limit points. 
We will denote a generic limit point by $\prob_{\bbeta}$ and its sampling measure by $\mu_{\bbeta}$.

Our goal is to express the free energy of the Hamiltonian
$$f(\bbeta):=\lim_{N\to\infty} \frac{1}{N}E \log Z_N(\bbeta)$$
as a functional of a limit point $\prob_{\bbeta}$.
For the mixed Hamiltonian with $p=1$ and even $p$'s, the existence of the limit was proved by Guerra and Toninelli~\cite{guerra_toninelli}.
We also know from the work of Guerra \cite{guerra_bound} and Talagrand \cite{tal_proof} that $f(\bbeta)$
is expressed by the celebrated Parisi formula.
We write it here in terms of Ruelle Probability Cascades (RPC's). The reader is referred to \cite{arguin_rpc} 
for the correspondence between this formulation and the one in terms of PDE's.
\begin{thm}[The Parisi Formula]
\label{thm: parisi formula}
For any $\bbeta$ satisfying \eqref{eqn: beta}, we have
$$
f(\bbeta)= \log 2 + \min_{\prob \text{ RPC} }\left\{ \mathcal{P}(\vec{\lambda}_{\bbeta}, \prob)-
\sum_p\frac{(p-1)\beta_p^2}{2}\Big(1-\E \mu ^{\otimes 2}(v\cdot v')^p\Big)\right\}
$$
where $\mathcal{P}$ is the Parisi functional \eqref{eqn:Parisi functional p} with 
$\psi(x)=\log\cosh x$,  and $(\vec{\lambda}_{\bbeta})_{p-1}=\sqrt{p}~ \beta_{p}$ for $p=1$
and even $p$'s and zero otherwise. The minimum is over the laws of the Ruelle Probability Cascades, a compact subset
 of the laws of ROSt's. The sampling measure of $\prob$ is denoted by $\mu$.
\end{thm}
It can be checked that for any ROSt that
is stochastically stable for the cavity fields $l_p$, so in particular for the RPC's, the second functional in the Parisi formula,
$$
\prob \mapsto \sum_p\frac{(p-1)\beta_p^2}{2}\Big(1-\E \mu ^{\otimes 2}(v\cdot v')^p\Big)
$$
is equal to the Parisi functional \eqref{eqn:Parisi functional p} for $\psi(x)=x$, $\lambda^2= \sum_p (p-1)\beta_p^2$,  and 
$c(x)=\sum_p\frac {(p-1)\beta_p^2}{\lambda^2}~x^p$.

The proof in \cite{tal_proof} of the Parisi Formula is not constructive. Our goal is to provide a constructive demonstration 
in the spirit of Aizenman, Sims and Starr \cite{ASS} and Theorem 2.4.19 of \cite{tal_book} by assuming that the Ultrametricity Conjecture holds.
This way of proceeding is based on general properties of ROSt's. 
Moreover, we emphasize that we will not use the ultrametric structure of the RPC explicitly, but simply that the
stochastic stability properties characterize the RPC family from the conjecture.

The first ingredient of the proof is an upper bound proved by Guerra \cite{guerra_bound}.
\begin{thm}[Guerra's bound]
For any RPC with law $\prob$ and sampling measure $\mu$, 
$$
f(\bbeta)\leq  \log 2 + \mathcal{P}(\vec{\lambda}_{\bbeta}, \prob)-
\sum_p\frac{(p-1)\beta_p^2}{2}\Big(1-\E \mu ^{\otimes 2}(v\cdot v')^p\Big)\ .
$$
\end{thm}
Guerra's proof is based on an interpolation between a cascade with a finite number of levels and
the Gibbs measure $\mathcal{G}_{\beta,N}$. It uses explicitly the tree structure of the cascades.
The proof of this bound can be also done using the approach of Aizenman, Sims and Starr \cite{ASS}.
In this case, the proof does not use the ultrametric structure of the cascade explicitly, but simply 
the stability of the cascades under the cavity mapping with function $\psi(x)=\log\cosh x$ and
$\psi(x)=x$. The reader is referred to \cite{ASS, arguin_rpc} for more details.

It remains to establish the matching lower bound.
We begin by a standard lemma, whose proof using Jensen's inequality will be helpful in this section.
\begin{lem}
\label{lem: continuity f}
The function $\bbeta\mapsto f(\bbeta)$ is convex and continuous
in the space of $\bbeta$ satisfying \eqref{eqn: beta}.
\end{lem}
\begin{proof}
We set $f_N(\bbeta):=\frac{1}{N}E \log Z_N(\bbeta)$. The fact that $f_N(\bbeta)$ is convex
follows by H\"older's inequality, and the convexity of $f(\bbeta)$ is proved since it is the pointwise
limit of $f_N(\bbeta)$. We now show that the family $(f_N(\bbeta))_N$ is equicontinuous in $\bbeta$.

First note that $f_N(\bbeta)$ is an increasing function of its parameters, since by Gaussian integration by part,
\begin{equation}
\label{eqn: deriv f}
D_p f_N(\beta)=\beta_p\Big(1-E \mathcal{G}_{\beta,N}^{\otimes 2}R(\sigma,\sigma')^p\Big)\geq 0\ ,
\end{equation}
where $D_p$ is the partial derivative in the $p$-th coordinate.
Let $\bbeta$ and $\bbeta'$ satisfy \eqref{eqn: beta}. Let us define $\bbeta^-$ with
$\beta^-_p:=\min\{\beta_p,\beta'_p\}$. By construction, $f_N(\bbeta^-)\leq f_N(\bbeta)$ and
$f_N(\bbeta^-)\leq f_N(\bbeta')$.
We show that
\begin{equation}
\label{eqn: f continuity}
f_N(\bbeta)-f_N(\bbeta^-)\leq \frac{1}{2}\sum_p (\beta_p-\beta_p^-)^2 \ ,
\end{equation}
and similarly for $f_N(\bbeta')$. We observe that
$$
f_N(\bbeta)-f_N(\bbeta^-)=\frac{1}{N} E E' \log \mathcal{G}_{\bbeta^-,N}(\exp H_N'(\sigma)),
$$
where $H_N'(\sigma)$ is defined as $H_N$, but is independent of $H_N$ and has parameters $\left(\sqrt{\beta_p^2-(\beta_p^-)^2}\right)_p$. 
We write $E'$ for the expectation on $H_N'$.
By Jensen's inequality with $\log$ and using the independence between $H_N'$ and $\mathcal{G}_{\bbeta^-}$, this is smaller than
$$
\frac{1}{N}E\Big[\log \mathcal{G}_{\bbeta^-,N} \Big( E' \exp H_N'(\sigma)\Big)\Big]= \frac{1}{2} \sum_{p} \big(\beta^2_p-(\beta_p^-)^2\big) \ ,
$$
since  $E' \exp H_N'(\sigma)=e^{\frac{N}{2}\big(\beta^2_p-(\beta_p^-)^2\big)}$. The proof is identical for $\bbeta'$.
This implies from the definition of $\bbeta^-$ that
$$
|f_N(\bbeta)-f_N(\bbeta')|\leq \frac{1}{2}\sum_p \big(\beta^2_p-(\beta'_p)^2\big)\leq \frac{1}{2}\sum_p 2^p \big(\beta^2_p-(\beta'_p)^2\big)\ .
$$
Since the continuity is uniform in $N$, the equicontinuity is proved. It follows that $f(\bbeta)$
is continuous in $\bbeta$ in the sense of the norm $\|\bbeta\|^2=\sum_p 2^p\beta_p^2$.
\end{proof}

It is well-known that a convex function on $\R$ is differentiable almost everywhere.
It turns out that a similar statement is true for a continuous convex functional on an infinite-dimensional space.
Indeed, it is a theorem of Mazur that a continuous convex functional on any open convex subset of a separable Banach space
has a dense set of points of differentiability (in the sense that the directional derivative
exists in all directions), see e.g. Theorem 1.20 in \cite{mazur}. Applying this to the 
Hilbert space of $\bbeta$ satisfying \eqref{eqn: beta}, we get
\begin{lem}
\label{lem: dense}
The set of parameters $\bbeta$ for which 
\begin{equation}
\label{eqn: derivative}
\lim_{t\to 0} \frac{f(\bbeta+t\bbeta')}{t} \text{ exists for any $\bbeta'$}
\end{equation}
is dense in the set of parameters satisfying \eqref{eqn: beta}.
\end{lem}
The main ingredient entering in the lower bound is the following result.
\begin{thm}
\label{thm: parisi 1}
Let $\bbeta$ be a point of differentiability of $f$.
There exists a subsequence of $(\prob_{\bbeta,N})$ that converges to $\prob_{\bbeta}$, the law of a ROSt with sampling
measure $\mu_{\bbeta}$, for which
\begin{equation}
\label{eqn: parisi}
f(\bbeta)\geq   \log 2 + \mathcal{P}(\vec{\lambda}_{\bbeta},\prob_{\bbeta} )-
\sum_p\frac{(p-1)\beta_p^2}{2}\Big(1-\E \mu_{\bbeta}^{\otimes 2}(v\cdot v')^p\Big)
\end{equation}
where $\mathcal{P}$ and $\vec{\lambda}_{\bbeta}$ are as in Theorem \ref{thm: parisi formula}.
\end{thm}

To prove the assertion, we shall need a standard result of convexity. We omit the proof.
\begin{lem}
\label{lem:convexity}
Let $f_N:\R^M\to\R$ be a sequence of differentiable convex functions that converges pointwise to $f$. Suppose $f$ is differentiable at $x\in\R^M$. 
Then for any sequence $x_N$ converging to $x$, $$\lim_{N\to\infty} \nabla f_N(x_N)= \lim_{N\to\infty}\nabla f_N(x)=\nabla f(x).$$
\end{lem}
\begin{proof}[Proof of Theorem \ref{thm: parisi 1}]
We show that 
\begin{equation}
\label{eqn: liminf}
\liminf_{N\to \infty}E\left[\log \frac{Z_{N+1}(\bbeta)}{Z_{N}(\bbeta)}\right]
\end{equation}
equals the right-hand side of \eqref{eqn: parisi} for $\prob_{\bbeta}$, a limit
point in the space of ROSt's of $(\prob_{\bbeta,N_k})_k$ where $(N_k)_k$ achieves the above $\liminf$.
The result then follows from the existence of the limiting free energy and the fact that
$$
\lim_{N\to\infty} f_N(\bbeta) =\lim_{N\to\infty}\frac{1}{N}\sum_{n=1}^NE\left[\log  \frac{Z_{n+1}(\bbeta)}{Z_{n}(\bbeta)}\right]\geq \liminf_{N\to \infty}E\left[\log \frac{Z_{N+1}(\bbeta)}{Z_{N}(\bbeta)}\right]\ .
$$

Taking $\bbeta^+:=(\beta_p^+)$ where $\beta_p^+=\beta_p\left(\frac{N+1}{N}\right)^{(p-1)/2}$, we can decompose the ratio in \eqref{eqn: liminf} as
\begin{equation}
\label{eqn:cavity step}
E\left[\log \frac{Z_{N+1}(\bbeta^+)}{Z_{N}(\bbeta)}\right]- E\left[\log \frac{Z_{N+1}(\bbeta^+)}{Z_{N+1}(\bbeta)}\right]\ .
\end{equation}
We handle the first term. A straightforward calculation using the representation \eqref{eqn: H_p} yields for $H_N^p
$ on $N+1$ spins,
$$
\begin{aligned}
\beta_p^+ H_{N+1}^p(\sigma, \sigma_{N+1})
&= \frac{\beta_p}{N^{(p-1)/2}} \sum_{k=0}^p \sigma_{N+1}^k \sum_{I_k}
  g_{I_k}\sigma_{j_1}...\sigma_{j_{p-k}}
\end{aligned}
$$
where the second sum is on $p$-tuplet $I_k=(i_1,...,i_{p})$ with exactly $k$ coordinates equal to $N+1$
and the other coordinates, denoted by $(j_1,...,j_{p-k})$, run from $1$ to $N$.
In particular, the covariance becomes
\begin{equation}
\label{eqn: cavity step 2}
\begin{aligned}
& E[ \beta_p^+ H_{N+1}^p(\sigma, \sigma_{N+1})\beta_p^+ H_{N+1}^p(\sigma', \sigma_{N+1}')]
= \beta^2_pE[H_N^p(\sigma)H_N^p(\sigma')]+
p\beta_p^2 \sigma_{N+1} \sigma'_{N+1} R(\sigma,\sigma')^{p-1}\\
&\hspace{8cm}+ \beta_p^2 \sum_{k=2}^p \binom{p}{k}  \sigma^k_{N+1} \sigma'^k_{N+1} \frac{R(\sigma,\sigma')^{p-k}}{N^{k-1}}
\ .
\end{aligned}
\end{equation}
We denote by $l_c(\sigma)$ the cavity field with covariance $c(R(\sigma,\sigma'))=\sum_p \frac{p\beta_p^2}{\lambda_{\bbeta}^2} R(\sigma,\sigma')^{p-1}$ where
$\lambda^2_{\bbeta}=\sum_p p\beta_p^2$, 
and by $\epsilon(\sigma,\sigma_{N+1})$
the field of covariance given by the sum over all $p$'s of the last term of \eqref{eqn: cavity step 2}.
Note that these fields are independent of each other and of $H_N(\sigma)$. 
Moreover, the variance of $\epsilon$ is smaller than
\begin{equation}
\label{eqn: epsilon}
\frac{1}{N}\sum_p\beta_p^2 \sum_{k=2}^p \binom{p}{k} \frac{1}{N^{k-2}}\leq \frac{1}{N}\sum_p \beta_p^2 2^p= \frac{C}{N}
\end{equation}
for some finite constant $C$ under the assumption \eqref{eqn: beta} on the decay of $\beta_p$.

The term $Z_{N+1}(\bbeta^+)$ can be rewritten using \eqref{eqn: cavity step 2} and we get
\begin{align*}
&E\left[\log \frac{Z_{N+1}(\bbeta^+)}{Z_{N}(\bbeta)}\right]\\
&= 
E \ E_l E_\epsilon \left[\log \G_{\beta,N}\left(\sum_{\sigma_{N+1}=\pm 1}\exp\{\lambda_{\bbeta}\sigma_{N+1} l_c(\sigma)+ \epsilon(\sigma,\sigma_{N+1})\}\right)\right],
\end{align*}
where $E$ is understood to be the expectation on $\G_{\beta,N}$ on the right-hand side, and $E_l$ and $E_\epsilon$ are the expectations on $l$ and $\epsilon$ respectively.
Using Jensen's inequality and the independence between $l$ and $\epsilon$ we get from \eqref{eqn: epsilon}
$$
\left| E\left[\log \frac{Z_{N+1}(\bbeta^+)}{Z_{N}(\bbeta)}\right]-
E E_l \left[\log \G_{\beta,N}\left(2\cosh \lambda_{\bbeta} l_c(\sigma)\right)\right]\right|
\leq \frac{C}{N}\ .
$$
Therefore the first term of \eqref{eqn:cavity step} is, up to a term of order $1/N$,
 $$
 \log 2+ \mathcal{P}(\vec{\lambda}_{\bbeta}, \prob_{\bbeta,N})\ ,
 $$ 
where  $ \prob_{\bbeta,N}$ is the law of the ROSt  constructed from $\G_{\bbeta,N}$. 
Let $\prob_{\bbeta}$ be a limit point of $(\prob_{\bbeta,N_k})$, where $(N_k)$ achieves the $\liminf$ of \eqref{eqn: liminf}.
The continuity of the Parisi functional proved in Proposition \ref{prop:continuity parisi} 
gives the first two terms of \eqref{eqn: parisi}. It remains to prove that along this subsequence
\begin{equation}
\label{eqn:cavity 2nd}
E\left[\log \frac{Z_{N_k+1}(\bbeta^+)}{Z_{N_k+1}(\bbeta)}\right]\to
\sum_p\frac{(p-1)\beta_p^2}{2}\Big(1-\E \mu_{\bbeta}^{\otimes 2}(v\cdot v')^p\Big)
\ .
\end{equation}
The idea is simply a Taylor expansion in $\bbeta$, but since we are dealing with infinitely many variables, we proceed with care.
We define $\bbeta_M^+$ with components $\beta_p^+$ for $p\leq M$ and $\beta_p$ for $p>M$.
Using Jensen's inequality and the summability of $\bbeta$, 
for any $\delta>0$ we can find $M(\delta)$ independently of $N$ such that for $M>M(\delta)$
$$
0\leq E\left[\log \frac{Z_{N+1}(\bbeta^+)}{Z_{N+1}(\bbeta_M^+)}\right]\leq \frac{N}{2}\sum_{p>M}({\beta_p^+}^2-\beta_p^2) < \delta\ .
$$
Fix such a $\delta$ and such a $M$.
By the mean-value theorem on $\R^M$, there exists $\overline{\bbeta}_M$ with $\beta_p\leq \bar{\beta}_p\leq \beta^+_p$
for $p\leq M$ and $\bar{\beta}_p=\beta_p$ for $p>M$ such that
$$
E\left[\log \frac{Z_{N+1}(\bbeta^+_M)}{Z_{N+1}(\bbeta)}\right]= \sum_{p\leq M} N(\beta_p^+ - \beta_p) D_p f_{N+1}(\overline{\bbeta}_M)\ .
$$
We know that $f_N(\bbeta)$ is a convex function of $\bbeta$.
Thus, by Lemma \ref{lem:convexity} and the assumption on differentiability of $f(\bbeta)$, 
$$
\begin{aligned}
\lim_{N\to\infty} \sum_{p\leq M} N(\beta_p^+ - \beta_p) D_p f_{N+1}(\overline{\bbeta}_M)
&=
\sum_{p\leq M} \frac{(p-1) \beta_p}{2} \lim_{N\to\infty} D_p f_N(\bbeta)\\
&= \sum_{p\leq M}\frac{(p-1)\beta_p^2}{2}\left(1-\E \mu_{\bbeta}^{\otimes 2} (v\cdot v')^p\right)\ ,
\end{aligned}
$$
where we used \eqref{eqn: deriv f} in the last equality as well as the continuity of the function $\prob\mapsto \E \mu^{\otimes 2} (v\cdot v')^p $ in the space of ROSt's. 
Since $\delta$ was arbitrary, the theorem is proved.
\end{proof}
The connection between the lower bound of Theorem \ref{thm: parisi 1} and the actual Parisi formula
with RPC's will be provided by the Ultrametricity Conjecture.
The idea is to use the lower bound of Theorem \ref{thm: parisi 1} for $\bbeta$ where $\beta_p>0$
for all $p$'s. The limit ROSt's $\prob_{\bbeta}$ for such a choice possesses strong invariance
properties by Corollary \ref{cor: diff->SS} if $\bbeta$ is a point of differentiability of $f$.
However, it is not guaranteed that $\bbeta$ is a point of differentiability.
Even worse, Mazur's theorem does not ensure that the points of differentiability are dense in the subset $\{\bbeta: \beta_p>0 ~ \forall p\}$
since it is not open. 
This apparent difficulty is however bypassed by the use of subsequences and by the continuity of the cavity mapping.
\begin{cor}
\label{cor: beta positive}
Let $\bbeta$ with $\beta_p>0$ for $p=1$ and even $p$'s. There exists a ROSt of law $\prob_{\bbeta}$
and sampling measure $\mu_{\bbeta}$
that is stochastically stable for the cavity fields $l_p$, $p=1$ and all $p$ even, and for which
$$
f(\bbeta)\geq   \log 2 + \mathcal{P}(\vec{\lambda}_{\bbeta},\prob_{\bbeta} )-
\sum_p\frac{(p-1)\beta_p^2}{2}\Big(1-\E \mu_{\bbeta}^{\otimes 2}(v\cdot v')^p\Big) \ .
$$
\end{cor}
\begin{proof}
By Lemma \ref{lem: dense}, there exists a sequence of points of differentiability $\bbeta_m$, $m\in\N$, that
converges to $\bbeta$. Hence the inequality \eqref{eqn: parisi} is satisfied for any $\bbeta_m$ for some ROSt's
of law $\prob_{\bbeta_m}$. Recall that $\prob_{\bbeta_m}$ is the limit law of the ROSt's constructed 
from the Gibbs measures $\G_{\bbeta_m,N}$. It follows from Corollary \ref{cor: diff->SS} that
$\prob_{\bbeta_m}$ is stochastically stable for the cavity fields $l_p$ for all $p\leq p(m)$, where
$p(m)\leq \infty$ is the greatest $p$ such that $(\bbeta_m)_p>0$. 
Since the space of ROSt's is compact, there exists a subsequence of $m$'s for which $\prob_{\bbeta_m}$ converges to
the law of some ROSt denoted $\prob_{\bbeta}$.
Moreover, since the cavity mapping is continuous, any limit point of a sequence of laws that are invariant
under a cavity mapping is also invariant. 
This implies that  $\prob_{\bbeta}$ is stochastically stable for the cavity fields $l_p$ for all relevant $p$'s
(since $p(m)$ goes to infinity with $m$ from the hypothesis that $\beta_p>0$ for all $p$).

By Lemma \ref{lem: continuity f}, $f(\bbeta_m)\to f(\bbeta)$. To prove the inequality  \eqref{eqn: parisi} for $\bbeta$, it remains to show that the
right-hand side also converges for the considered subsequence of $\bbeta_m$.
This is clear for the second term since $\bbeta_m\to\bbeta$ and $\prob_{\bbeta_m}\to \prob_{\bbeta}$. 
For the second term, one has
$$
\begin{aligned}
\left|\mathcal{P}(\vec{\lambda}_{\bbeta_m},\prob_{\bbeta_m} )-\mathcal{P}(\vec{\lambda}_{\bbeta},\prob_{\bbeta} )\right|
&\leq \left|\mathcal{P}(\vec{\lambda}_{\bbeta_m},\prob_{\bbeta_m} )-\mathcal{P}(\vec{\lambda}_{\bbeta},\prob_{\bbeta_m} )\right|
+ \left|\mathcal{P}(\vec{\lambda}_{\bbeta},\prob_{\bbeta_m} )-\mathcal{P}(\vec{\lambda}_{\bbeta},\prob_{\bbeta} )\right|\\
&\leq \max_{\prob \text{ ROSt}}\left|\mathcal{P}(\vec{\lambda}_{\bbeta_m},\prob )-\mathcal{P}(\vec{\lambda}_{\bbeta},\prob )\right|
+ \left|\mathcal{P}(\vec{\lambda}_{\bbeta},\prob_{\bbeta_m} )-\mathcal{P}(\vec{\lambda}_{\bbeta},\prob_{\bbeta} )\right|\ .
\end{aligned}
$$
The right-hand side goes to zero along the subsequence of $\bbeta_m$ by Proposition \ref{prop:continuity parisi}
because $\prob_{\bbeta_m}\to \prob_{\bbeta}$ and $\| \vec{\lambda}_{\bbeta_m}-\vec{\lambda}_{\bbeta}\|_2=\sum_pp((\bbeta_m)^2_p-\beta_p^2)\leq \| \bbeta_m-\bbeta\|$.
\end{proof}

\begin{proof}[Proof of Theorem \ref{thm: parisi formula} using Conjecture \ref{conj:ultra}]
By Conjecture \ref{conj:ultra}, the only stochastically stable ROSt's for infinitely many $p$'s are convex combinations of RPC's.
Therefore, by Corollary \ref{cor: beta positive} and Guerra's bound,
the Parisi formula is established at any point $\bbeta$ with $\beta_p>0$ for all $p$. 
(Note that the functionals that are minimized in the Parisi formula are linear functions in the law of the ROSt's. 
In particular, the minimum can be restricted to RPC's as opposed to convex combinations of RPC's.)
It remains to show that the Parisi formula holds at all $\bbeta$, including those $\bbeta$ with $\beta_p=0$ for some $p$. 
Since $f$ is continuous in the space of parameters satisfying \eqref{eqn: beta}, it suffices to show that the functional
 \begin{equation}
 \label{eqn: Lambda}
 \bbeta \mapsto \min_{\prob \text{ RPC} }\left\{ \mathcal{P}(\vec{\lambda}_{\bbeta}, \prob)-
\sum_p\frac{(p-1)\beta_p^2}{2}\Big(1-\E \mu ^{\otimes 2}(v\cdot v')^p\Big)\right\}
 \end{equation}
 is continuous under the norm $\|\bbeta\|^2=\sum_p 2^p \beta_p^2$.
We write for convenience $\Lambda(\bbeta, \prob)$ for the functional in the argument of the minimum.
By Proposition \ref{prop:continuity parisi}, $\Lambda(\bbeta, \cdot)$ 
is a continuous functional on ROSt's. We denote by $\prob_{\bbeta}$ the RPC where the minimum of $\Lambda(\bbeta,\cdot)$ is attained.
Since we minimize $\Lambda$ over ROSt's, we have trivially
\begin{equation}
\label{eqn: estimate continuity}
 \Lambda(\bbeta,\prob_{\bbeta}) - \Lambda(\bbeta',\prob_{\bbeta}) 
 \leq \min_{\prob \text{ RPC}}\Lambda(\bbeta,\prob) - \min_{\prob \text{ RPC}}\Lambda(\bbeta',\prob)\leq
  \Lambda(\bbeta,\prob_{\bbeta'}) - \Lambda(\bbeta',\prob_{\bbeta'})\ .
\end{equation}
On the other hand, if $\bbeta$ is close enough to $\bbeta'$, we have for some $K>0$,
$$
\max_{\prob \text{ ROSt}} \left|\Lambda(\bbeta, \prob)-\Lambda(\bbeta', \prob)\right|\leq
K\|\bbeta-\bbeta' \| \ .
$$
This is because, by the second part of Proposition \ref{prop:continuity parisi},
$$
\max_{\prob \text{ ROSt}} \left|\mathcal{P}(\vec{\lambda}_{\bbeta}, \prob)-\mathcal{P}(\vec{\lambda}_{\bbeta'}, \prob)\right|\leq K \left\{\sum_{p} p(\beta_p-\beta'_p)^2\right\}^{1/2}
\leq  K \left\{\sum_p 2^p(\beta_p-\beta'_p)^2\right\}^{1/2}
  \ ,
$$
and the same holds by inspection for the second functional entering in the definition of $\Lambda$.
The continuity of \eqref{eqn: Lambda}
in the space of $\bbeta$ satisfying \eqref{eqn: beta} follows from \eqref{eqn: estimate continuity} and the above estimate.
\end{proof}


\section{Properties of Stochastically Stable ROSt's}
\label{section:properties}
Some particular properties of stochastically stable ROSt's were 
essential in the proof of the ultrametricity conjecture in \cite{argaiz} and \cite{panchenko_GG},
where strong assumptions on the ROSt were needed (e.g., finiteness of the support for the law of the entries).
Here we show similar properties in a general settting which may be useful for the proof of the conjecture in the general case.

A direct consequence of the continuity of the cavity mapping is the compactness of the set of stochastically stable ROSt's.
\begin{cor}
\label{cor:compact}
The set of laws of ROSt's that are invariant under a (possibly infinite) 
family of cavity mappings of the form \eqref{eqn:mapping2} is convex and compact. 
In particular, the set of laws of the stochastically stable ROSt's is convex and compact.
\end{cor}
\begin{proof}
Convexity is clear. The compactness of the set of laws of ROSt's that are invariant under the mapping \eqref{eqn:mapping2} for
a given $\psi$ and a given positive definite form $c$ is a consequence of continuity of the mapping (Theorem \ref{thm:continuity}).
Compactness also holds for sets of ROSt's that are invariant under an infinite number of mappings (e.g., for all parameters $\lambda>0$) since in a compact space
an arbitrary intersection of compact sets is compact.
\end{proof}

Corollary \ref{cor:compact} yields a standard ergodic decomposition for 
the different sets of stochastically stable ROSt's by Choquet's theorem:
the law of a ROSt in such a convex set can be written as 
the barycenter of some probability measure on the extreme points of the set. 
This is useful since many functions of interest on the sampling measure of a ROSt turn out to be constant
when the law is an extreme point.
This is helpful when proving properties of the sampling measure of stochastically stable ROSt's.

\begin{df}
\label{df: invariant fct}
We say a function $h:\M(\B)\to \R$ is an {\it invariant function of the cavity mapping} $\Phi_{\psi,l_c}$ defined in \eqref{eqn:mapping2} if:
1) it is Borel measurable; 2) for any isometry $T$ of $\B$, $h(\mu)=h(T^{-1}\mu)$; 3) $h(\mu)=h(\Phi_{\psi,l_c}\mu)$ $P_l$-almost surely. 
\end{df}

\begin{lem}
\label{lem:dim}
Let $h:\M(\B)\to \R$ be an invariant function of a (possibly infinite) family of cavity mapping $(\Phi_{\psi, l_c})$, where the collection is over different $\psi$'s and $c$'s.
Let $\prob$ be the law of a ROSt with sampling measure $\mu$.
If $\prob$ is an extreme point of the set of laws that are invariant under the family $(\Phi_{\psi, l_c})$, 
then the random variable $h(\mu)$ is constant $\prob$-a.s.
\end{lem}

\begin{proof}
If $h$ is not a constant, we have the non-trivial decomposition of $\prob$
\begin{equation}
\label{eqn:decomp}
\prob= \int_{\R}\prob(\ \cdot \ | \ h(\mu)=y)\ \prob h^{-1}(dy)\ .
\end{equation}
On the other hand, $\prob(\ \cdot \ | h(\mu))$ is invariant under the same family of mapping  $(\Phi_{\psi, l_c})$ for $\prob$-almost all $h(\mu)$.
 Indeed, for any Borel measurable function $g$ on $\M(\B)$, invariant under isometry, any Borel measurable function $F:\R\to\R$,
and every $\Phi_{\psi, l_c}$ in the collection,
$$
\E[g(\mu)F\big(h(\mu)\big) ]= \E E_l [g(\Phi_{\psi, l_c} \mu)F\big(h(\Phi_{\psi, l_c} \mu)\big) ]=  \E E_l [g(\Phi_{\psi, l_c} \mu)F\big(h(\mu)\big) ]\ ,
$$
where the first equality follows from the invariance of $\prob$ and the second, from the invariance of $h$.
Thus \eqref{eqn:decomp} contradicts the assumption that $\prob$ is an extreme point and the claim follows.
\end{proof}
By inspection of the form of the cavity mapping \eqref{eqn:mapping2}, 
we expect that the support of a measure $\mu$ is not affected by the cavity mapping,
since the mapping modifies only the weight of the vectors and not their relative position.
We make this idea precise in the following lemma.
\begin{lem}
\label{lem: equivalence}
Let $\mu\in\M(\B)$ and $\Phi_{\psi,l_c}$ as in \eqref{eqn:mapping2}. Then $\Phi_{\psi,l_c} \mu$ is equivalent to $\mu$ $P_l$-almost surely, that is:
for any non-negative measurable function $f$ on $\B$,
$$
\int_{\B}f(v)\ \mu(dv)=0 \ \text{ if and only if }\ \int_\B f(v) \Phi_{\psi,l_c}\mu(dv)=0 \text{ $P_l$-a.s.}
$$
\end{lem}
\begin{proof}
Note first that
$$
\int_\B f(v) \Phi_{\psi,l_c}\mu(dv)=\int_\B  \frac{f(v) ~ E_z[e^{\psi(l_c( v)+z\sqrt{1-c(\| v\|^2)})}]}{\mu(E_z[e^{\psi(l_c( v)+z\sqrt{1-c(\| v\|^2)})}])}\mu(dv)
$$
hence the left-hand side is zero if and only if $\int_\B f(v) E_z[e^{\psi(l_c( v)+z\sqrt{1-c(\| v\|^2)})}] \mu(dv)=0$.
Now if $\int_{\B}f(v)\ \mu(dv)=0$, then we must have that $f(v)=0$ $\mu$-a.e., which implies 
$$\int_\B f(v) E_z[e^{\psi(l_c( v)+z\sqrt{1-c(\| v\|^2)})}] \mu(dv)=0\ .$$
On the other hand, we note that $E_z[e^{\psi(l_c( v)+z\sqrt{1-c(\| v\|^2)})}] >0$ $P_l$-a.s. and $\mu$-a.e. 
This is because $l_c( v)+z\sqrt{1-c(\| v\|^2)}$ is finite $P_l\times P_z$-a.s. and $\mu$-a.e., since by Fubini, 
$$E_l E_z\int_\B |l_c( v)+z\sqrt{1-c(\| v\|^2)}|^2\ \mu(dv)= 1\ .$$
Therefore if $\int_\B f(v) E_z[e^{\psi(l_c( v)+z\sqrt{1-c(\| v\|^2)})}] \mu(dv)=0$, then $f(v)=0$ $\mu$-a.e. also.
\end{proof}
Functions on $\M(\B)$ that only depend on the support of the measure typically turn out to be invariant functions of the mapping in the sense of Definition \ref{df: invariant fct}. 
This simple fact can be used to investigate the support of the sampling measures of stochastically stable ROSt's.
Here we look at two examples below: the dimension of the support of the measure and 
whether or not the support lies on a sphere. 
We will need three invariant functions of the cavity mapping.

The first two are
\begin{equation}
\label{eqn: r max min}
\begin{aligned}
r_{\min}(\mu)&:= \inf \Big\{ 0\leq r\leq 1 : \  \mu\{v\in\B: \|v\| \leq r\}>0\Big\} \ ;\\
r_{\max}(\mu)&:=\sup \Big\{ 0\leq r\leq 1: \ \mu\{v\in\B:\|v\|\geq r\}>0\Big\}\ .
\end{aligned}
\end{equation}
Essentially, $r_{\min}$ represents the inner radius of the support of $\mu$ and $r_{\max}$, the outer radius.
The requirements of Definition \ref{df: invariant fct} are fulfilled by these two functions for any cavity mapping $\Phi_{\psi,l_c}$.
Indeed, both functions are measurable since
\begin{equation}
\label{eqn: r}
\begin{aligned}
\{\mu\in\M(\B):r_{\min}(\mu)<s\}&=\bigcup_{s'\in\Q,\ 0\leq s'<s}\Big\{\mu\in\M(\B): \mu\{\|v\|\leq s'\}>0\Big\}\\
\{\mu\in\M(\B):r_{\max}(\mu)>s\}&=\bigcup_{s'\in\Q,\ s<s'\leq 1}\Big\{\mu\in\M(\B): \mu\{\|v\|\geq s'\}>0\Big\}\ ,
\end{aligned}
\end{equation}
and the functions $\mu\mapsto \mu\{\|v\|\leq s'\}$, $\mu\mapsto \mu\{\|v\|\geq s'\}$ are Borel measurable for the weak-* topology on $\M(\B)$.
They are invariant under an isometry $T$ of $\B$ since the norm of a vector is. 
The fact that $r_{\min}(\mu)=r_{\min}(\Phi_{\psi, l_c}\mu)$ and $r_{\max}(\mu)=r_{\max}(\Phi_{\psi, l_c}\mu)$ is direct from Lemma \ref{lem: equivalence}, since the latter implies $P_l$-a.s.
$$
\begin{aligned}
\mu\{v\in\B: \|v\| \leq s\}>0 &\Longleftrightarrow \Phi_{\psi, l_c} \mu\{v\in\B: \|v\| \leq s\}>0 \\
\mu\{v\in\B: \|v\| \geq s\}>0 &\Longleftrightarrow \Phi_{\psi, l_c} \mu\{v\in\B: \|v\| \geq s\}>0 \ .
\end{aligned}
$$

The third invariant function we need is the dimension of the smallest subspace that contains the support.
One way to make this precise is as follows \cite{aizenman_priv}. Let $\mu\in\M(\B)$. Define (weakly) the covariance operator $\hat{\mathcal{C}}_\mu$ 
on $\hilbert$ as
$$
\text{for all $v',v''\in\hilbert$, } ~ v'\cdot (\hat{\mathcal{C}}_\mu v''):= \int_{\B}\mu(dv)\ (v'\cdot v)(v\cdot v'') \ .
$$
It is easily checked that $\hat{\mathcal{C}}_\mu$ is a self-adjoint, trace-class linear operator on $\hilbert$.
(It is trace-class, since for a standard basis $\{e_i\}$ of $\hilbert$, 
$Tr \ \hat{\mathcal{C}}_\mu=\sum_{i} e_i \cdot (\hat{\mathcal{C}}_\mu e_i)=\int_{\B}\mu(dv)\|v\|^2\leq 1$.)
In particular, it is compact, thus admits a basis of orthonormal eigenvectors for $\hilbert$. We write $\{\lambda_i(\mu)\}$ 
for the non-zero eigenvalues. 
We write $\hilbert_\mu$ for the eigenspace corresponding to the non-zero eigenvalues.
We have $\hilbert=\hilbert_\mu\oplus \hilbert_\mu^{\perp}$. 
It is easily verified that $\mu(\hilbert_\mu^{\perp}\cap\B)=0$.
We define
$$
\dim \mu:=\#\{i\in\N: \lambda_i(\mu)>0\}=\dim \hilbert_\mu\ .
$$
We claim that the function $\mu\mapsto \dim \mu$ is an invariant of the cavity mapping $\Phi_{\psi,l_c}$ for any 
$\psi$ and $c$.
Standard arguments show that it is measurable with respect to the Borel $\sigma$-algebra of $\M(\B)$.
It is readily checked that eigenvectors of $\mu$ are mapped to eigenvectors of $\mu T^{-1}$ for any isometry $T$ of $\B$,
therefore $\dim \mu =\dim \mu T^{-1}$.
It remains to show that $\dim \mu=\dim \Phi_{\psi,l_c}\mu$ $P_l$-a.s.
We have $v\in \hilbert_{\mu}^{\perp}\cap\B$ if and only if
$$
v\cdot(\hat{\mathcal{C}}_\mu v )= \int_{\B}\mu(dv)\ (v,v)^2=0\ .
$$
By Lemma \ref{lem: equivalence}, this is true if and only if
$$
\int_{\B}\Phi_{\psi,l_c}\mu(dv)\ (v,v)^2=v\cdot(\hat{\mathcal{C}}_{\Phi_{\psi,l_c}\mu}\ v) =0, 
$$
which in turn is equivalent to $v\in \hilbert_{\Phi_{\psi,l_c}\mu}^{\perp}\cap\B$. The claim follows.

The following result on the dimension of stochastically stable ROSt is a generalization of a result for competing particle systems which states that no system with a finite number of particles can be quasi-stationary \cite{ruzmaiz}. It is rather surprising that it holds for ROSt whose a law is invariant under a single cavity mapping.
\begin{thmdim}
Let $\prob$ be the law of a ROSt that is invariant under the cavity mapping $\Phi_{\lambda,l}$
for a given $\lambda>0$. 
Let $\mu$ be its sampling measure.
Then 
$$
\prob(\dim \mu=1 \text{ or } \dim \mu=\infty)=1 \ .
$$
In other words, $\mu$ is supported on a single vector or on an infinite-dimensional subset of $\B$.
\end{thmdim}

\begin{proof}
For simplicity, we set $\lambda=1$.
By relying on the Choquet decomposition, it suffices to show the assertion holds when $\prob$ is an extreme point of the set of laws invariant under $\Phi_{1,l}$.
Since $\mu\mapsto \dim \mu$  is an invariant function of the mapping, Lemma \ref{lem:dim} implies that it is constant with $\prob$-probability one. 
Suppose $1<\dim \mu<\infty$. 
We first prove that this implies that the sampling measure $\mu$ is supported on a sphere $\prob$-almost surely, that is, there exists a deterministic $r\geq 0$ 
such that $\mu\{v\in\B:\|v\|=r\}=1$ for almost all $\mu$.
Indeed, since $\mu\mapsto r_{\min}(\mu)$ defined in \eqref{eqn: r} is also an invariant function, 
it is constant $\prob$-a.s.
This constant will henceforth be written simply as $r_{\min}$. 
For any $r>r_{\min}$ and $\epsilon>0$, consider the ratio
\begin{equation}
\label{eqn: ratio}
\frac{\mu\{v\in\B: \|v\|> r + \epsilon \} }{\mu\{v\in\B: \|v\|\leq r  \} }\ .
\end{equation}
We claim that the ratio is zero for any $r>r_{\min}$ and $\epsilon>0$.
If so, we are done, since it implies that $\mu\{v\in\B: \|v\|>r\}=0$ for any $r>r_{\min}$, hence that $\mu$ is supported on the sphere
of radius $r_{\min}$.
By Lemma \ref{lem: composition}, $\prob$ is invariant under the mappings $\Phi_{\sqrt{T},l}$
for any $T\in\N$.
Therefore the ratio \eqref{eqn: ratio} has the same law as
$$
\frac{\int_{\{\|v\|> r + \epsilon\}} e^{\sqrt{T}l(v) -\frac{T}{2}\|v\|^2} \mu(dv) }
{\int_{\{\|v\|\leq r\}} e^{\sqrt{T}l(v) -\frac{T}{2}\|v\|^2} \mu(dv) }\ .
$$
Since $1<\dim \mu<\infty$, $l$ can be seen as a Gaussian vector in a finite-dimensional space. We then have the bound 
$-\|l\|\leq l(v) \leq \|l\|$, since $\|v\|\leq 1$. 
A straightforward application of that bound and simple estimates shows that the ratio is smaller or equal to
$$
e^{2\sqrt{T} \|l\|-\frac{T}{2}(2\epsilon r +\epsilon^2)}\frac{\mu\{v\in\B: \|v\|> r + \epsilon \} }{\mu\{v\in\B: \|v\|\leq r  \} }\ .
$$
This proves the claim, since the above goes to zero as $T\to\infty$ for $P_l$-almost all $l$.

It remains to show that if $\mu$ has support in a finite-dimensional sphere, then it must be supported on a single vector.
The idea is that $\mu$ must be supported on the vector where the cavity field is maximal.
In a finite-dimensional space, $v\mapsto l(v)$ is a continuous function for $P_l$-almost all $l$.
Since the sphere is compact, it achieves its maximum, say $m_l$. Moreover, by independence of the $l(v)$'s in orthogonal directions, 
the maximizer, say $v_l^*$, is unique.
Let $C_l(\epsilon)$ be a cap neighborhood of the maximizer $v_l^*$ on the sphere such that $l(v)>m_l -\epsilon$ 
for $v \in C_l(\epsilon)$ 
and $l(v)\leq m_l -\epsilon$ for $v \in C_l(\epsilon)^c$. If the measure $\mu$ is not supported on a
single vector, there exists an $\epsilon$ such that $\int_{C_l(\epsilon)^c} e^{\sqrt{T} l(v)} \mu(dv)>0$. We then have
$$
\frac{\int_{C_l(\epsilon/2)} e^{\sqrt{T} l(v)}\mu(dv) }{\int_{C_l(\epsilon)^c} e^{\sqrt{T} l(v)} \mu(dv) }
\geq e^{\sqrt{T}\epsilon/2}\frac{\mu\{C_l(\epsilon/2)\}}{\mu\{C_l(\epsilon)^c\}}\ .
$$ 
This shows that, for any $\epsilon$, the limit of the left-hand side as $T\to\infty$ exists and is infinite 
for $\prob$-almost all $\mu$ and $P_l$-almost all $l$.
But by invariance, the ratio on the left-hand side has the same law under $\prob\times P_l$ as for a single application of $\Phi_{1,l}$, i.e., 
$$
\frac{\int_{C_l(\epsilon/2)} e^{ l(v)}\mu(dv) }{\int_{C_l(\epsilon)^c} e^{ l(v)} \mu(dv)}\ .
$$
This ratio is infinite for all $\epsilon$ if and only if
$$
(\Phi_{1,l} \mu)\{v_l^*\}=1 \text{ \  $\prob\times P_l$-a.s.}
$$
We deduce by the invariance of the ROSt that $\mu$ must be supported on a single vector. 
\end{proof}

One might expect that other properties of stochastically stable ROSt's might be obtained 
in a similar way leading to a better understanding of the hypotheses of the Ultrametricity Conjecture.
If stochastic stability is assumed for an infinite number of cavity fields, we can prove more.

\begin{thmsphere}
Let $\prob$ be the law of a ROSt that is invariant under the cavity mapping $\Phi_{\lambda,l_p}$
for a given $\lambda>0$ and for infinitely many $p\in\N$.
Let $\mu$ be its sampling measure and suppose that $\mu^{\otimes 2}\{(v,v')\in\B^2:  v\cdot v'=r^2_{\max}(\mu)\}>0$ 
on a set of positive $\prob$-probability.
Then,
$$
\mu\{v\in\B: \|v\|=r_{\max}(\mu)\}=1\ \text{ $\prob$-a.s.}
$$
In other words,  $\mu$ is almost surely supported on a sphere.
\end{thmsphere}

\begin{proof}
We set $\lambda=1$ for simplicity. In the proof, all $p$'s refer to the infinite collection
of $p$'s for which invariance holds.
By relying again on the Choquet decomposition, 
it suffices to show the assertion holds when $\prob$ is an extreme point of the set of laws invariant under $\Phi_{1,l_p}$
for all $p$'s.
Consider $r_{\min}$ and $r_{\max}$ defined in  \eqref{eqn: r}.
Since they are invariant functions of the cavity mappings $\Phi_{1,l_p}$ for all $p$,
they are constant $\prob$-a.s. 
We suppose that $r_{\min}<r_{\max}$ and show this leads to a contradiction.

Pick $\delta$ such that $r_{\min}<\delta<r_{\max}$. Let $B_{\delta}=\{v\in\B:\|v\|<\delta\}$. 
Note that since $r_{\min}<r_{\max}$, we must have that $0<\mu(B_\delta)<1$ for $\prob$-almost all $\mu$.

By the invariance properties of the ROSt and Lemma \ref{lem: composition}, the following identity holds
for all $T\in\N$ and all $p$'s
\begin{equation}
\label{eqn: inv p}
\E\mu(B_\delta)=
\E E_l\left[\frac{\int_{B_\delta} e^{\sqrt{T} \ l_p(v)-\frac{T}{2}\|v\|^{2p}} \mu(dv) }{\mu\left(e^{\sqrt{T}\ l_p(v)-\frac{T}{2}\|v\|^{2p}}\right)}\right]\ .
\end{equation}
Now pick a sequence $(T(p))_p$ such that $T(p)\delta^{2p}\to 0$. 
(The choice of the sequence will be refined later.)
Fatou's lemma then implies
\begin{equation}
\label{eqn: fatou}
\E\mu(B_\delta)\geq \E \liminf_{p\to\infty} E_l\left[\frac{\int_{B_\delta} e^{\sqrt{T(p)} \ l_p(v)-\frac{T(p)}{2}\|v\|^{2p}} \mu(dv) }{\mu\left(e^{\sqrt{T(p)}\ l_p(v)-\frac{T(p)}{2}\|v\|^{2p}}\right)}\right]
\end{equation}
We now show that for $\prob$-almost all $\mu$,
\begin{equation}
\label{eqn: limit p}
 \liminf_{p\to\infty} E_l\left[\frac{\int_{B_\delta} e^{\sqrt{T(p)} \ l_p(v)-\frac{T(p)}{2}\|v\|^{2p}} \mu(dv) }{\mu\left(e^{\sqrt{T(p)}\ l_p(v)-\frac{T(p)}{2}\|v\|^{2p}}\right)}\right]
 =
\liminf_{p\to\infty} E_l\left[\frac{ \mu(B_\delta) }{\mu\left(e^{\sqrt{T(p)}\ l_p(v)-\frac{T(p)}{2}\|v\|^{2p}}\right)}\right]\ .
\end{equation}
Taking the difference of the two sides and applying Cauchy-Schwarz inequality, one gets
$$
\left\{ E_l\left(\int_{B_\delta} e^{\sqrt{T(p)} \ l_p(v)-\frac{1}{2}T(p)\|v\|^{2p}} \mu(dv)-\mu(B_\delta)\right)^2 \right\}^{1/2}
\left\{ E_l\left[\frac{1}{\left(\mu\left(e^{\sqrt{T(p)}\ l_p(v)-\frac{T(p)}{2}\|v\|^{2p}}\right)\right)^2}\right] \right\}^{1/2}\ .
$$
The term inside the first bracket can be evaluated by developing the square and integrating over $E_l$ to get
$$
 \int_{B_\delta\times B_\delta} e^{T(p) (v\cdot v')^{p}} \mu(dv)\mu(dv') -(\mu(B_\delta))^2\ ,
$$
which converges to $0$ as $p\to\infty$ by the choice of $T(p)$.
It remains to show that the term in the second bracket is bounded uniformly in $p$. By introducing $\mu(B_\delta)$, it is equal to
$$
\frac{1}{(\mu(B_\delta))^2}E_l\left[\left(\frac{\mu(B_\delta)}{\int_{B_\delta}e^{\sqrt{T(p)}\ l_p(v)-\frac{T(p)}{2}\|v\|^{2p}}\mu(dv)}\right)^2\right]
$$
(Note that we used the fact that $\mu(B_\delta)>0$ here.) The term in the expectation can be written in terms of the conditional measure $\mu(~|B_\delta)$.
Hence using Jensen's inequality, one gets the upper bound
$$
\frac{1}{(\mu(B_\delta))^2}E_l \mu\left( e^{-2\sqrt{T(p)}\ l_p(v)+T(p)\|v\|^{2p}} \big| B_\delta \right)\ ,
$$ 
which is seen to be bounded by $e^{3T(p)\delta^{2p}}$ by integrating with $E_l$.

Putting \eqref{eqn: fatou} and \eqref{eqn: limit p} together yields
$$
\E\mu(B_\delta)\geq \E\left[  \mu(B_\delta)  \liminf_{p\to\infty} E_l\left[\frac{1}{\mu\left(e^{\sqrt{T(p)}\ l_p(v)-\frac{T(p)}{2}\|v\|^{2p}}\right)}\right]\right]\ .
$$
And since
$$
E_l\left[\frac{1}{\mu\left(e^{\sqrt{T(p)}\ l_p(v)-\frac{T(p)}{2}\|v\|^{2p}}\right)}\right]\geq  \left(E_l\mu\left(e^{\sqrt{T(p)}\ l_p(v)-\frac{T(p)}{2}\|v\|^{2p}}\right)\right)^{-1}=1, 
$$
and $\mu(B_\delta)>0$ $\prob$-a.s., this implies that
\begin{equation}
\label{eqn: contradiction}
\liminf_{p\to\infty} E_l\left[\frac{1}{\mu\left(e^{\sqrt{T(p)}\ l_p(v)-\frac{T(p)}{2}\|v\|^{2p}}\right)}\right]=1 \ \text{ $\prob$-a.s.}
\end{equation}
By hypothesis, there must exist an $\alpha>0$ such that the set of $\mu$'s with
$\mu^{\otimes 2}\{(v,v'): v\cdot v'=r_{\max}^2\}>\alpha$ has positive probability.
We show that the equality \eqref{eqn: contradiction} cannot hold on this set.
Fix a $\mu$ in this set.
Write for simplicity $m_p:=\mu\left(e^{\sqrt{T(p)}\ l_p(v)-\frac{T(p)}{2}\|v\|^{2p}}\right)$.
In view of the fact that the function $y\mapsto y^{-1}$ is strictly convex,
a contradiction would be reached if we can construct along a subsequence of 
$(m_p)_p$ a random variable $m$ that is fluctuating but for which $E_lm=E_lm^{-1}=1$.
It follows directly from the value of the expectations and Markov's inequality that for any $\epsilon>0$, there exists $0<\eta<1$ such that
$$
P_l\{m_p\in [\eta,\eta^{-1}]\} \geq 1-\epsilon\ .
$$
This implies that the sequence of random variables $(m_p)$ is tight. Moreover, proceeding as above, 
$$
E_lm_p^{-2}\leq \frac{1}{(\mu(B_\delta))^2}E_l\mu\left( e^{-2\sqrt{T(p)}\ l_p(v)+T(p)\|v\|^{2p}} \big| B_\delta \right)<\infty \text{ , uniformly in $p$.}
$$
Therefore there exists a subsequence of $p$'s such that $m_p\to m$ in law and $E_l m_p^{-1}\to E_lm^{-1}=1$
by \eqref{eqn: contradiction} and by the uniform boundedness of the above moment. It remains to show that $E_lm=1$ and that the variance of $m$ is non-zero so that $m$ fluctuates.
Integrating over $E_l$ one has the identity
$$
E_l m_p^2= \mu^{\otimes 2}\left(e^{T(p)(v\cdot v')^p}\right)\ .
$$
To prove both claims, it suffices to show that the limit of $E_l m_p^2$ along the subsequence is finite (ensuring the limit of $E_lm_p$ is $E_lm$) and strictly greater than one 
(ensuring the variance is non-zero). This will be true if we choose the integers $T(p)$ such that $T(p)r_{\max}^{2p}\to \log \frac{2}{\alpha}$. 
(Note that $T(p)\delta^{2p}\to 0$ for any $\delta<r_{\max}$ for this choice.)
In this case, the following hold
$$
\begin{aligned}
E_l m_p^2=&\mu^{\otimes 2}\left(e^{T(p)(v\cdot v')^p}\right)\geq \alpha e^{T(p)r_{\max}^{2p}} \to 2 >1\\
E_l m_p^2=&\mu^{\otimes 2}\left(e^{T(p)(v\cdot v')^p}\right)\leq e^{T(p)r_{\max}^{2p}} \to \frac{2}{\alpha}< \infty\ ,
\end{aligned}
$$
leading to the contradiction.
\end{proof}

 We remark that stochastic stability alone is not sufficient to imply that the sampling measure of a ROSt be supported on a sphere.
This is demonstrated by the following example. Take $p=(p_i)_{i\in\N}$ a Poisson-Dirichlet variable 
with parameter $x_1$, and $p'=(p'_j)_{j\in\N}$ a Poisson-Dirichlet variable with parameter $x_2\neq x_1$. 
Pick $q_1\neq q_2$ strictly positive such that $q_1+q_2=1$ and 
\begin{equation}\label{eq:ex support}
(1-x_1)q_1=(1-x_2)q_2\ .
\end{equation}
It suffices to consider the ROSt with sampling measure $\mu=\sum_{i}p_i\ \delta_{u_i}+\sum_{j}p'_j\ \delta_{v_j}$ where $(u_i)$ are orthogonal vectors of norm $\sqrt{q_1}$ and $(v_j)$ are orthogonal vectors with norm $\sqrt{q_2}$
also orthogonal to the $u$'s. 
A Poisson-Dirichlet variable with parameter $x$ has the property that, for any positive random variable $W$
such that $E W^x<\infty$ and for a iid sequence $(W_i)$ with the law of $W$ (see e.g. \cite{argaiz}),
$$(p_iW_i)\distrib (E[W^x]^{1/x}p_i)\ .$$
Therefore taking $W_i=e^{l(u_i)}$ and $W_j=e^{l(v_j)}$, we have for the variables defined above and for the cavity field $l$:
$$
\begin{aligned}
(p_ie^{l(u_i)-\frac{q_1}{2}})\distrib e^{\frac{-(1-x_1) q_1}{2}}(p_i) \hspace{0.2cm} \text{, and } \hspace{0.2cm}(p'_je^{l(v_j)-\frac{q_2}{2}})\distrib e^{\frac{-(1-x_2) q_2}{2}}(p'_j)\ .
\end{aligned}
$$
The relation \eqref{eq:ex support} ensures that the constant is the same in both cases, hence vanishes after normalization of the weights.
This shows that the constructed $\mu$ is stochastically stable

The above setup also provides a counterexample to the Ultrametricity Conjecture \ref{conj:ultra}
if stochastic stability for a single cavity field alone holds.
This is demonstrated by the following example coming from the limit of two uncoupled REM's as studied in \cite{lp_kistler}. 
It suffices to consider the ROSt with sampling measure
$\mu=\sum_{i,j}p_i p'_j \ \delta_{\frac{1}{\sqrt{2}}(u_i+v_j)}$. It is easily checked that this ROSt is stochastically stable from the invariance property of Poisson-Dirichlet variables 
and the fact that $l(\frac{1}{\sqrt{2}}(u_i+v_j))=\frac{1}{\sqrt{2}}l(u_i)+\frac{1}{\sqrt{2}}l(v_j)$ 
(the restriction \eqref{eq:ex support} is not needed here). 
On the other hand, the support of its sampling measure is not ultrametric as it can be easily checked from the 
following triplet of vectors in the support of $\mu$: $\frac{1}{\sqrt{2}}(u_1+v_1)$, $\frac{1}{\sqrt{2}}(u_2+v_2)$ and $\frac{1}{\sqrt{2}}(u_1+v_2)$.

\end{document}